\newread\epsffilein    
\newif\ifepsffileok    
\newif\ifepsfbbfound   
\newif\ifepsfverbose   
\newdimen\epsfxsize    
\newdimen\epsfysize    
\newdimen\epsftsize    
\newdimen\epsfrsize    
\newdimen\epsftmp      
\newdimen\pspoints     
\pspoints=1bp          
\epsfxsize=0pt         
\epsfysize=0pt         
\def\epsfbox#1{\global\def\epsfllx{72}\global\def\epsflly{72}%
   \global\def\epsfurx{540}\global\def\epsfury{720}%
   \def\lbracket{[}\def\testit{#1}\ifx\testit\lbracket
   \let\next=\epsfgetlitbb\else\let\next=\epsfnormal\fi\next{#1}}%
\def\epsfgetlitbb#1#2 #3 #4 #5]#6{\epsfgrab #2 #3 #4 #5 .\\%
   \epsfsetgraph{#6}}%
\def\epsfnormal#1{\epsfgetbb{#1}\epsfsetgraph{#1}}%
\def\epsfgetbb#1{%
%
%
\openin\epsffilein=#1
\ifeof\epsffilein\errmessage{I couldn't open #1, will ignore it}\else
%
%
   {\epsffileoktrue \chardef\other=12
    \def\do##1{\catcode`##1=\other}\dospecials \catcode`\ =10
    \loop
       \read\epsffilein to \epsffileline
       \ifeof\epsffilein\epsffileokfalse\else
%
%
          \expandafter\epsfaux\epsffileline:. \\%
       \fi
   \ifepsffileok\repeat
   \ifepsfbbfound\else
    \ifepsfverbose\message{No bounding box comment in #1; using defaults}\fi\fi
   }\closein\epsffilein\fi}%
%
%
\def\epsfsetgraph#1{%
   \epsfrsize=\epsfury\pspoints
   \advance\epsfrsize by-\epsflly\pspoints
   \epsftsize=\epsfurx\pspoints
   \advance\epsftsize by-\epsfllx\pspoints
%
%
   \epsfxsize\epsfsize\epsftsize\epsfrsize
   \ifnum\epsfxsize=0 \ifnum\epsfysize=0
      \epsfxsize=\epsftsize \epsfysize=\epsfrsize
%
%
     \else\epsftmp=\epsftsize \divide\epsftmp\epsfrsize
       \epsfxsize=\epsfysize \multiply\epsfxsize\epsftmp
       \multiply\epsftmp\epsfrsize \advance\epsftsize-\epsftmp
       \epsftmp=\epsfysize
       \loop \advance\epsftsize\epsftsize \divide\epsftmp 2
       \ifnum\epsftmp>0
          \ifnum\epsftsize<\epsfrsize\else
             \advance\epsftsize-\epsfrsize \advance\epsfxsize\epsftmp \fi
       \repeat
     \fi
   \else\epsftmp=\epsfrsize \divide\epsftmp\epsftsize
     \epsfysize=\epsfxsize \multiply\epsfysize\epsftmp   
     \multiply\epsftmp\epsftsize \advance\epsfrsize-\epsftmp
     \epsftmp=\epsfxsize
     \loop \advance\epsfrsize\epsfrsize \divide\epsftmp 2
     \ifnum\epsftmp>0
        \ifnum\epsfrsize<\epsftsize\else
           \advance\epsfrsize-\epsftsize \advance\epsfysize\epsftmp \fi
     \repeat     
   \fi
%
%
   \ifepsfverbose\message{#1: width=\the\epsfxsize, height=\the\epsfysize}\fi
   \epsftmp=10\epsfxsize \divide\epsftmp\pspoints
   \vbox to\epsfysize{\vfil\hbox to\epsfxsize{%
      \includegraphics{#1}%
      \hfil}}%
\epsfxsize=0pt\epsfysize=0pt}%

%
%
{\catcode`\%=12 \global\let\epsfpercent=
%
%
\long\def\epsfaux#1#2:#3\\{\ifx#1\epsfpercent
   \def\testit{#2}\ifx\testit\epsfbblit
      \epsfgrab #3 . . . \\%
      \epsffileokfalse
      \global\epsfbbfoundtrue
   \fi\else\ifx#1\par\else\epsffileokfalse\fi\fi}%
%
%
\def\epsfgrab #1 #2 #3 #4 #5\\{%
   \global\def\epsfllx{#1}\ifx\epsfllx\empty
      \epsfgrab #2 #3 #4 #5 .\\\else
   \global\def\epsflly{#2}%
   \global\def\epsfurx{#3}\global\def\epsfury{#4}\fi}%
%
%
\def\epsfsize#1#2{\epsfxsize}

\catcode`!=11 
 
  

\def\PiC{P\kern-.12em\lower.5ex\hbox{I}\kern-.075emC}
\def\PiCTeX{\PiC\kern-.11em\TeX}

\def\!ifnextchar#1#2#3{%
  \let\!testchar=#1%
  \def\!first{#2}%
  \def\!second{#3}%
  \futurelet\!nextchar\!testnext}
\def\!testnext{%
  \ifx \!nextchar \!spacetoken 
    \let\!next=\!skipspacetestagain
  \else
    \ifx \!nextchar \!testchar
      \let\!next=\!first
    \else 
      \let\!next=\!second 
    \fi 
  \fi
  \!next}
\def\\{\!skipspacetestagain} 
  \expandafter\def\\ {\futurelet\!nextchar\!testnext} 
\def\\{\let\!spacetoken= } \\  

\def\!tfor#1:=#2\do#3{%
  \edef\!fortemp{#2}%
  \ifx\!fortemp\!empty 
    \else
    \!tforloop#2\!nil\!nil\!!#1{#3}%
  \fi}
\def\!tforloop#1#2\!!#3#4{%
  \def#3{#1}%
  \ifx #3\!nnil
    \let\!nextwhile=\!fornoop
  \else
    #4\relax
    \let\!nextwhile=\!tforloop
  \fi 
  \!nextwhile#2\!!#3{#4}}

\def\!etfor#1:=#2\do#3{%
  \def\!!tfor{\!tfor#1:=}%
  \edef\!!!tfor{#2}%
  \expandafter\!!tfor\!!!tfor\do{#3}}

\def\!cfor#1:=#2\do#3{%
  \edef\!fortemp{#2}%
  \ifx\!fortemp\!empty 
  \else
    \!cforloop#2,\!nil,\!nil\!!#1{#3}%
  \fi}
\def\!cforloop#1,#2\!!#3#4{%
  \def#3{#1}%
  \ifx #3\!nnil
    \let\!nextwhile=\!fornoop 
  \else
    #4\relax
    \let\!nextwhile=\!cforloop
  \fi
  \!nextwhile#2\!!#3{#4}}

\def\!ecfor#1:=#2\do#3{%
  \def\!!cfor{\!cfor#1:=}%
  \edef\!!!cfor{#2}%
  \expandafter\!!cfor\!!!cfor\do{#3}}

\def\!empty{}
\def\!nnil{\!nil}
\def\!fornoop#1\!!#2#3{}

\def\!ifempty#1#2#3{%
  \edef\!emptyarg{#1}%
  \ifx\!emptyarg\!empty
    #2%
  \else
    #3%
  \fi}
 
\def\!getnext#1\from#2{%
  \expandafter\!gnext#2\!#1#2}%
\def\!gnext\\#1#2\!#3#4{%
  \def#3{#1}%
  \def#4{#2\\{#1}}%
  \ignorespaces}

%
\def\!getnextvalueof#1\from#2{%
  \expandafter\!gnextv#2\!#1#2}%
\def\!gnextv\\#1#2\!#3#4{%
  #3=#1%
  \def#4{#2\\{#1}}%
  \ignorespaces}

\def\!copylist#1\to#2{%
  \expandafter\!!copylist#1\!#2}
\def\!!copylist#1\!#2{%
  \def#2{#1}\ignorespaces}

\def\!wlet#1=#2{%
  \let#1=#2 
  \wlog{\string#1=\string#2}}
 
\def\!listaddon#1#2{%
  \expandafter\!!listaddon#2\!{#1}#2}
\def\!!listaddon#1\!#2#3{%
  \def#3{#1\\#2}}
 

\def\!rightappend#1\withCS#2\to#3{\expandafter\!!rightappend#3\!#2{#1}#3}
\def\!!rightappend#1\!#2#3#4{\def#4{#1#2{#3}}}

\def\!leftappend#1\withCS#2\to#3{\expandafter\!!leftappend#3\!#2{#1}#3}
\def\!!leftappend#1\!#2#3#4{\def#4{#2{#3}#1}}

\def\!lop#1\to#2{\expandafter\!!lop#1\!#1#2}
\def\!!lop\\#1#2\!#3#4{\def#4{#1}\def#3{#2}}



\def\!loop#1\repeat{\def\!body{#1}\!iterate}
\def\!iterate{\!body\let\!next=\!iterate\else\let\!next=\relax\fi\!next}
 
\def\!!loop#1\repeat{\def\!!body{#1}\!!iterate}
\def\!!iterate{\!!body\let\!!next=\!!iterate\else\let\!!next=\relax\fi\!!next}
 
\def\!removept#1#2{\edef#2{\expandafter\!!removePT\the#1}}
{\catcode`p=12 \catcode`t=12 \gdef\!!removePT#1pt{#1}}

\def\placevalueinpts of <#1> in #2 {%
  \!removept{#1}{#2}}
 
\def\!mlap#1{\hbox to 0pt{\hss#1\hss}}
\def\!vmlap#1{\vbox to 0pt{\vss#1\vss}}
 
\def\!not#1{%
  #1\relax
    \!switchfalse
  \else
    \!switchtrue
  \fi
  \if!switch
  \ignorespaces}


 

\let\!!!wlog=\wlog              
\def\wlog#1{}    

\newdimen\headingtoplotskip     
\newdimen\linethickness         
\newdimen\longticklength        
\newdimen\plotsymbolspacing     
\newdimen\shortticklength       
\newdimen\stackleading          
\newdimen\tickstovaluesleading  
\newdimen\totalarclength        
\newdimen\valuestolabelleading  

\newbox\!boxA                   
\newbox\!boxB                   
\newbox\!picbox                 
\newbox\!plotsymbol             
\newbox\!putobject              
\newbox\!shadesymbol            

\newcount\!countA               
\newcount\!countB               
\newcount\!countC               
\newcount\!countD               
\newcount\!countE               
\newcount\!countF               
\newcount\!countG               
\newcount\!fiftypt              
\newcount\!intervalno           
\newcount\!npoints              
\newcount\!nsegments            
\newcount\!ntemp                
\newcount\!parity               
\newcount\!scalefactor          
\newcount\!tfs                  
\newcount\!tickcase             

\newdimen\!Xleft                
\newdimen\!Xright               
\newdimen\!Xsave                
\newdimen\!Ybot                 
\newdimen\!Ysave                
\newdimen\!Ytop                 
\newdimen\!angle                
\newdimen\!arclength            
\newdimen\!areabloc             
\newdimen\!arealloc             
\newdimen\!arearloc             
\newdimen\!areatloc             
\newdimen\!bshrinkage           
\newdimen\!checkbot             
\newdimen\!checkleft            
\newdimen\!checkright           
\newdimen\!checktop             
\newdimen\!dimenA               
\newdimen\!dimenB               
\newdimen\!dimenC               
\newdimen\!dimenD               
\newdimen\!dimenE               
\newdimen\!dimenF               
\newdimen\!dimenG               
\newdimen\!dimenH               
\newdimen\!dimenI               
\newdimen\!distacross           
\newdimen\!downlength           
\newdimen\!dp                   
\newdimen\!dshade               
\newdimen\!dxpos                
\newdimen\!dxprime              
\newdimen\!dypos                
\newdimen\!dyprime              
\newdimen\!ht                   
\newdimen\!leaderlength         
\newdimen\!lshrinkage           
\newdimen\!midarclength         
\newdimen\!offset               
\newdimen\!plotheadingoffset    
\newdimen\!plotsymbolxshift     
\newdimen\!plotsymbolyshift     
\newdimen\!plotxorigin          
\newdimen\!plotyorigin          
\newdimen\!rootten              
\newdimen\!rshrinkage           
\newdimen\!shadesymbolxshift    
\newdimen\!shadesymbolyshift    
\newdimen\!tenAa                
\newdimen\!tenAc                
\newdimen\!tenAe                
\newdimen\!tshrinkage           
\newdimen\!uplength             
\newdimen\!wd                   
\newdimen\!wmax                 
\newdimen\!wmin                 
\newdimen\!xB                   
\newdimen\!xC                   
\newdimen\!xE                   
\newdimen\!xM                   
\newdimen\!xS                   
\newdimen\!xaxislength          
\newdimen\!xdiff                
\newdimen\!xleft                
\newdimen\!xloc                 
\newdimen\!xorigin              
\newdimen\!xpivot               
\newdimen\!xpos                 
\newdimen\!xprime               
\newdimen\!xright               
\newdimen\!xshade               
\newdimen\!xshift               
\newdimen\!xtemp                
\newdimen\!xunit                
\newdimen\!xxE                  
\newdimen\!xxM                  
\newdimen\!xxS                  
\newdimen\!xxloc                
\newdimen\!yB                   
\newdimen\!yC                   
\newdimen\!yE                   
\newdimen\!yM                   
\newdimen\!yS                   
\newdimen\!yaxislength          
\newdimen\!ybot                 
\newdimen\!ydiff                
\newdimen\!yloc                 
\newdimen\!yorigin              
\newdimen\!ypivot               
\newdimen\!ypos                 
\newdimen\!yprime               
\newdimen\!yshade               
\newdimen\!yshift               
\newdimen\!ytemp                
\newdimen\!ytop                 
\newdimen\!yunit                
\newdimen\!yyE                  
\newdimen\!yyM                  
\newdimen\!yyS                  
\newdimen\!yyloc                
\newdimen\!zpt                  

\newif\if!axisvisible           
\newif\if!gridlinestoo          
\newif\if!keepPO                
\newif\if!placeaxislabel        
\newif\if!switch                
\newif\if!xswitch               

\newtoks\!axisLaBeL             
\newtoks\!keywordtoks           

\newwrite\!replotfile           

\newhelp\!keywordhelp{The keyword mentioned in the error message in unknown. 
Replace NEW KEYWORD in the indicated response by the keyword that 
should have been specified.}    

\!wlet\!!origin=\!xM                   
\!wlet\!!unit=\!uplength               
\!wlet\!Lresiduallength=\!dimenG       
\!wlet\!Rresiduallength=\!dimenF       
\!wlet\!axisLength=\!distacross        
\!wlet\!axisend=\!ydiff                
\!wlet\!axisstart=\!xdiff              
\!wlet\!axisxlevel=\!arclength         
\!wlet\!axisylevel=\!downlength        
\!wlet\!beta=\!dimenE                  
\!wlet\!gamma=\!dimenF                 
\!wlet\!shadexorigin=\!plotxorigin     
\!wlet\!shadeyorigin=\!plotyorigin     
\!wlet\!ticklength=\!xS                
\!wlet\!ticklocation=\!xE              
\!wlet\!ticklocationincr=\!yE          
\!wlet\!tickwidth=\!yS                 
\!wlet\!totalleaderlength=\!dimenE     
\!wlet\!xone=\!xprime                  
\!wlet\!xtwo=\!dxprime                 
\!wlet\!ySsave=\!yM                    
\!wlet\!ybB=\!yB                       
\!wlet\!ybC=\!yC                       
\!wlet\!ybE=\!yE                       
\!wlet\!ybM=\!yM                       
\!wlet\!ybS=\!yS                       
\!wlet\!ybpos=\!yyloc                  
\!wlet\!yone=\!yprime                  
\!wlet\!ytB=\!xB                       
\!wlet\!ytC=\!xC                       
\!wlet\!ytE=\!downlength               
\!wlet\!ytM=\!arclength                
\!wlet\!ytS=\!distacross               
\!wlet\!ytpos=\!xxloc                  
\!wlet\!ytwo=\!dyprime                 

\!zpt=0pt                              
\!xunit=1pt
\!yunit=1pt
\!arearloc=\!xunit
\!areatloc=\!yunit
\!dshade=5pt
\!leaderlength=24in
\!tfs=256                              
\!wmax=5.3pt                           
\!wmin=2.7pt                           
\!xaxislength=\!xunit
\!xpivot=\!zpt
\!yaxislength=\!yunit 
\!ypivot=\!zpt
\plotsymbolspacing=.4pt
  \!dimenA=50pt \!fiftypt=\!dimenA     

\!rootten=3.162278pt                   
\!tenAa=8.690286pt                     
\!tenAc=2.773839pt                     
\!tenAe=2.543275pt                     

\def\!cosrotationangle{1}      
\def\!sinrotationangle{0}      
\def\!xpivotcoord{0}           
\def\!xref{0}                  
\def\!xshadesave{0}            
\def\!ypivotcoord{0}           
\def\!yref{0}                  
\def\!yshadesave{0}            
\def\!zero{0}                  

\let\wlog=\!!!wlog
%
  
\def\normalgraphs{%
  \longticklength=.4\baselineskip
  \shortticklength=.25\baselineskip
  \tickstovaluesleading=.25\baselineskip
  \valuestolabelleading=.8\baselineskip
  \linethickness=.4pt
  \stackleading=.17\baselineskip
  \headingtoplotskip=1.5\baselineskip
  \visibleaxes
  \ticksout
  \nogridlines
  \unloggedticks}
%
\def\setplotarea x from #1 to #2, y from #3 to #4 {%
  \!arealloc=\!M{#1}\!xunit \advance \!arealloc -\!xorigin
  \!areabloc=\!M{#3}\!yunit \advance \!areabloc -\!yorigin
  \!arearloc=\!M{#2}\!xunit \advance \!arearloc -\!xorigin
  \!areatloc=\!M{#4}\!yunit \advance \!areatloc -\!yorigin
  \!initinboundscheck
  \!xaxislength=\!arearloc  \advance\!xaxislength -\!arealloc
  \!yaxislength=\!areatloc  \advance\!yaxislength -\!areabloc
  \!plotheadingoffset=\!zpt
  \!dimenput {{\setbox0=\hbox{}\wd0=\!xaxislength\ht0=\!yaxislength\box0}}
     [bl] (\!arealloc,\!areabloc)}
%
\def\visibleaxes{%
  \def\!axisvisibility{\!axisvisibletrue}}

%

\def\!fixkeyword#1{%
  \errhelp=\!keywordhelp
  \errmessage{Unrecognized keyword `#1': \the\!keywordtoks{NEW KEYWORD}'}}

\!keywordtoks={enter `i\fixkeyword}

\def\fixkeyword#1{%
  \!nextkeyword#1 }


\def\axis {%
  \def\!nextkeyword##1 {%
    \expandafter\ifx\csname !axis##1\endcsname \relax
      \def\!next{\!fixkeyword{##1}}%
    \else
      \def\!next{\csname !axis##1\endcsname}%
    \fi
    \!next}%
  \!offset=\!zpt
  \!axisvisibility
  \!placeaxislabelfalse
  \!nextkeyword}

\def\!axisbottom{%
  \!axisylevel=\!areabloc
  \def\!tickxsign{0}%
  \def\!tickysign{-}%
  \def\!axissetup{\!axisxsetup}%
  \def\!axislabeltbrl{t}%
  \!nextkeyword}

\def\!axistop{%
  \!axisylevel=\!areatloc
  \def\!tickxsign{0}%
  \def\!tickysign{+}%
  \def\!axissetup{\!axisxsetup}%
  \def\!axislabeltbrl{b}%
  \!nextkeyword}

\def\!axisleft{%
  \!axisxlevel=\!arealloc
  \def\!tickxsign{-}%
  \def\!tickysign{0}%
  \def\!axissetup{\!axisysetup}%
  \def\!axislabeltbrl{r}%
  \!nextkeyword}

\def\!axisright{%
  \!axisxlevel=\!arearloc
  \def\!tickxsign{+}%
  \def\!tickysign{0}%
  \def\!axissetup{\!axisysetup}%
  \def\!axislabeltbrl{l}%
  \!nextkeyword}

\def\!axisshiftedto#1=#2 {%
  \if 0\!tickxsign
    \!axisylevel=\!M{#2}\!yunit
    \advance\!axisylevel -\!yorigin
  \else
    \!axisxlevel=\!M{#2}\!xunit
    \advance\!axisxlevel -\!xorigin
  \fi
  \!nextkeyword}

\def\!axisvisible{%
  \!axisvisibletrue  
  \!nextkeyword}

\def\!axisinvisible{%
  \!axisvisiblefalse
  \!nextkeyword}

\def\!axislabel#1 {%
  \!axisLaBeL={#1}%
  \!placeaxislabeltrue
  \!nextkeyword}

\expandafter\def\csname !axis/\endcsname{%
  \!axissetup 
  \if!placeaxislabel
    \!placeaxislabel
  \fi
  \if +\!tickysign 
    \!dimenA=\!axisylevel
    \advance\!dimenA \!offset 
    \advance\!dimenA -\!areatloc 
    \ifdim \!dimenA>\!plotheadingoffset
      \!plotheadingoffset=\!dimenA 
    \fi
  \fi}

\def\grid #1 #2 {%
  \!countA=#1\advance\!countA 1
  \axis bottom invisible ticks length <\!zpt> andacross quantity {\!countA} /
  \!countA=#2\advance\!countA 1
  \axis left   invisible ticks length <\!zpt> andacross quantity {\!countA} / }

\def\plotheading#1 {%
  \advance\!plotheadingoffset \headingtoplotskip
  \!dimenput {#1} [B] <.5\!xaxislength,\!plotheadingoffset>
    (\!arealloc,\!areatloc)}

\def\!axisxsetup{%
  \!axisxlevel=\!arealloc
  \!axisstart=\!arealloc
  \!axisend=\!arearloc
  \!axisLength=\!xaxislength
  \!!origin=\!xorigin
  \!!unit=\!xunit
  \!xswitchtrue
  \if!axisvisible 
    \!makeaxis
  \fi}

\def\!axisysetup{%
  \!axisylevel=\!areabloc
  \!axisstart=\!areabloc
  \!axisend=\!areatloc
  \!axisLength=\!yaxislength
  \!!origin=\!yorigin
  \!!unit=\!yunit
  \!xswitchfalse
  \if!axisvisible
    \!makeaxis
  \fi}

\def\!makeaxis{%
  \setbox\!boxA=\hbox{
    \beginpicture
      \!setdimenmode
      \setcoordinatesystem point at {\!zpt} {\!zpt}   
      \putrule from {\!zpt} {\!zpt} to
        {\!tickysign\!tickysign\!axisLength} 
        {\!tickxsign\!tickxsign\!axisLength}
    \endpicturesave <\!Xsave,\!Ysave>}%
    \wd\!boxA=\!zpt
    \!placetick\!axisstart}

\def\!placeaxislabel{%
  \advance\!offset \valuestolabelleading
  \if!xswitch
    \!dimenput {\the\!axisLaBeL} [\!axislabeltbrl]
      <.5\!axisLength,\!tickysign\!offset> (\!axisxlevel,\!axisylevel)
    \advance\!offset \!dp  
    \advance\!offset \!ht  
  \else
    \!dimenput {\the\!axisLaBeL} [\!axislabeltbrl]
      <\!tickxsign\!offset,.5\!axisLength> (\!axisxlevel,\!axisylevel)
  \fi
  \!axisLaBeL={}}

%


\def\arrow <#1> [#2,#3]{%
  \!ifnextchar<{\!arrow{#1}{#2}{#3}}{\!arrow{#1}{#2}{#3}<\!zpt,\!zpt> }}

\def\!arrow#1#2#3<#4,#5> from #6 #7 to #8 #9 {%
%
  \!xloc=\!M{#8}\!xunit   
  \!yloc=\!M{#9}\!yunit
  \!dxpos=\!xloc  \!dimenA=\!M{#6}\!xunit  \advance \!dxpos -\!dimenA
  \!dypos=\!yloc  \!dimenA=\!M{#7}\!yunit  \advance \!dypos -\!dimenA
  \let\!MAH=\!M
  \!setdimenmode
  \!xshift=#4\relax  \!yshift=#5\relax
  \!reverserotateonly\!xshift\!yshift
  \advance\!xshift\!xloc  \advance\!yshift\!yloc
%
  \!xS=-\!dxpos  \advance\!xS\!xshift
  \!yS=-\!dypos  \advance\!yS\!yshift
  \!start (\!xS,\!yS)
  \!ljoin (\!xshift,\!yshift)
%
  \!Pythag\!dxpos\!dypos\!arclength
  \!divide\!dxpos\!arclength\!dxpos  
  \!dxpos=32\!dxpos  \!removept\!dxpos\!!cos
  \!divide\!dypos\!arclength\!dypos  
  \!dypos=32\!dypos  \!removept\!dypos\!!sin
%
  \!halfhead{#1}{#2}{#3}
  \!halfhead{#1}{-#2}{-#3}
  \let\!M=\!MAH
  \ignorespaces}
%
  \def\!halfhead#1#2#3{%
    \!dimenC=-#1%
    \divide \!dimenC 2 
    \!dimenD=#2\!dimenC
    \!rotate(\!dimenC,\!dimenD)by(\!!cos,\!!sin)to(\!xM,\!yM)
    \!dimenC=-#1
    \!dimenD=#3\!dimenC
    \!dimenD=.5\!dimenD
    \!rotate(\!dimenC,\!dimenD)by(\!!cos,\!!sin)to(\!xE,\!yE)
    \!start (\!xshift,\!yshift)
    \advance\!xM\!xshift  \advance\!yM\!yshift
    \advance\!xE\!xshift  \advance\!yE\!yshift
    \!qjoin (\!xM,\!yM) (\!xE,\!yE) 
    \ignorespaces}

\def\betweenarrows #1#2 from #3 #4 to #5 #6 {%
  \!xloc=\!M{#3}\!xunit  \!xxloc=\!M{#5}\!xunit%
  \!yloc=\!M{#4}\!yunit  \!yyloc=\!M{#6}\!yunit%
  \!dxpos=\!xxloc  \advance\!dxpos by -\!xloc
  \!dypos=\!yyloc  \advance\!dypos by -\!yloc
  \advance\!xloc .5\!dxpos
  \advance\!yloc .5\!dypos
  \let\!MBA=\!M
  \!setdimenmode
  \ifdim\!dypos=\!zpt
    \ifdim\!dxpos<\!zpt \!dxpos=-\!dxpos \fi
    \put {\!lrarrows{\!dxpos}{#1}}#2{} at {\!xloc} {\!yloc}
  \else
    \ifdim\!dxpos=\!zpt
      \ifdim\!dypos<\!zpt \!dypos=-\!dypos \fi
      \put {\!udarrows{\!dypos}{#1}}#2{} at {\!xloc} {\!yloc}
    \fi
  \fi
  \let\!M=\!MBA
  \ignorespaces}

\def\!lrarrows#1#2{
  {\setbox\!boxA=\hbox{$\mkern-2mu\mathord-\mkern-2mu$}%
   \setbox\!boxB=\hbox{$\leftarrow$}\!dimenE=\ht\!boxB
   \setbox\!boxB=\hbox{}\ht\!boxB=2\!dimenE
   \hbox to #1{$\mathord\leftarrow\mkern-6mu
     \cleaders\copy\!boxA\hfil
     \mkern-6mu\mathord-$%
     \kern.4em $\vcenter{\box\!boxB}$$\vcenter{\hbox{#2}}$\kern.4em
     $\mathord-\mkern-6mu
     \cleaders\copy\!boxA\hfil
     \mkern-6mu\mathord\rightarrow$}}}

\def\!udarrows#1#2{
  {\setbox\!boxB=\hbox{#2}%
   \setbox\!boxA=\hbox to \wd\!boxB{\hss$\vert$\hss}%
   \!dimenE=\ht\!boxA \advance\!dimenE \dp\!boxA \divide\!dimenE 2
   \vbox to #1{\offinterlineskip
      \vskip .05556\!dimenE
      \hbox to \wd\!boxB{\hss$\mkern.4mu\uparrow$\hss}\vskip-\!dimenE
      \cleaders\copy\!boxA\vfil
      \vskip-\!dimenE\copy\!boxA
      \vskip\!dimenE\copy\!boxB\vskip.4em
      \copy\!boxA\vskip-\!dimenE
      \cleaders\copy\!boxA\vfil
      \vskip-\!dimenE \hbox to \wd\!boxB{\hss$\mkern.4mu\downarrow$\hss}
      \vskip .05556\!dimenE}}}

%

\def\putbar#1breadth <#2> from #3 #4 to #5 #6 {%
  \!xloc=\!M{#3}\!xunit  \!xxloc=\!M{#5}\!xunit%
  \!yloc=\!M{#4}\!yunit  \!yyloc=\!M{#6}\!yunit%
  \!dypos=\!yyloc  \advance\!dypos by -\!yloc
  \!dimenI=#2  
  \ifdim \!dimenI=\!zpt 
    \putrule#1from {#3} {#4} to {#5} {#6} 
  \else 
    \let\!MBar=\!M
    \!setdimenmode 
    \divide\!dimenI 2
    \ifdim \!dypos=\!zpt             
      \advance \!yloc -\!dimenI 
      \advance \!yyloc \!dimenI
    \else
      \advance \!xloc -\!dimenI 
      \advance \!xxloc \!dimenI
    \fi
    \putrectangle#1corners at {\!xloc} {\!yloc} and {\!xxloc} {\!yyloc}
    \let\!M=\!MBar 
  \fi
  \ignorespaces}

\def\setbars#1breadth <#2> baseline at #3 = #4 {%
  \edef\!barshift{#1}%
  \edef\!barbreadth{#2}%
  \edef\!barorientation{#3}%
  \edef\!barbaseline{#4}%
  \def\!bardobaselabel{\!bardoendlabel}%
  \def\!bardoendlabel{\!barfinish}%
  \let\!drawcurve=\!barcurve
  \!setbars}
\def\!setbars{%
  \futurelet\!nextchar\!!setbars}
\def\!!setbars{%
  \if b\!nextchar
    \def\!!!setbars{\!setbarsbget}%
  \else 
    \if e\!nextchar
      \def\!!!setbars{\!setbarseget}%
    \else
      \def\!!!setbars{\relax}%
    \fi
  \fi
  \!!!setbars}
\def\!setbarsbget baselabels (#1) {%
  \def\!barbaselabelorientation{#1}%
  \def\!bardobaselabel{\!!bardobaselabel}%
  \!setbars}
\def\!setbarseget endlabels (#1) {%
  \edef\!barendlabelorientation{#1}%
  \def\!bardoendlabel{\!!bardoendlabel}%
  \!setbars}

\def\!barcurve #1 #2 {%
  \if y\!barorientation
    \def\!basexarg{#1}%
    \def\!baseyarg{\!barbaseline}%
  \else
    \def\!basexarg{\!barbaseline}%
    \def\!baseyarg{#2}%
  \fi
  \expandafter\putbar\!barshift breadth <\!barbreadth> from {\!basexarg}
    {\!baseyarg} to {#1} {#2}
  \def\!endxarg{#1}%
  \def\!endyarg{#2}%
  \!bardobaselabel}

\def\!!bardobaselabel "#1" {%
  \put {#1}\!barbaselabelorientation{} at {\!basexarg} {\!baseyarg}
  \!bardoendlabel}
 
\def\!!bardoendlabel "#1" {%
  \put {#1}\!barendlabelorientation{} at {\!endxarg} {\!endyarg}
  \!barfinish}

\def\!barfinish{%
  \!ifnextchar/{\!finish}{\!barcurve}}

%
%
%
\def\putrectangle{%
  \!ifnextchar<{\!putrectangle}{\!putrectangle<\!zpt,\!zpt> }}
\def\!putrectangle<#1,#2> corners at #3 #4 and #5 #6 {%
%
  \!xone=\!M{#3}\!xunit  \!xtwo=\!M{#5}\!xunit%
  \!yone=\!M{#4}\!yunit  \!ytwo=\!M{#6}\!yunit%
  \ifdim \!xtwo<\!xone
    \!dimenI=\!xone  \!xone=\!xtwo  \!xtwo=\!dimenI
  \fi
  \ifdim \!ytwo<\!yone
    \!dimenI=\!yone  \!yone=\!ytwo  \!ytwo=\!dimenI
  \fi
  \!dimenI=#1\relax  \advance\!xone\!dimenI  \advance\!xtwo\!dimenI
  \!dimenI=#2\relax  \advance\!yone\!dimenI  \advance\!ytwo\!dimenI
  \let\!MRect=\!M
  \!setdimenmode
%
  \!shaderectangle
%
  \!dimenI=.5\linethickness
  \advance \!xone  -\!dimenI
  \advance \!xtwo   \!dimenI
  \putrule from {\!xone} {\!yone} to {\!xtwo} {\!yone} 
  \putrule from {\!xone} {\!ytwo} to {\!xtwo} {\!ytwo} 
%
  \advance \!xone   \!dimenI
  \advance \!xtwo  -\!dimenI%
  \advance \!yone  -\!dimenI
  \advance \!ytwo   \!dimenI
  \putrule from {\!xone} {\!yone} to {\!xone} {\!ytwo} 
  \putrule from {\!xtwo} {\!yone} to {\!xtwo} {\!ytwo} 
  \let\!M=\!MRect
  \ignorespaces}
 

\def\shaderectanglesoff{%
  \def\!shaderectangle{}%
  \ignorespaces}

\shaderectanglesoff
 
\def\!!shaderectangle{%
  \!dimenA=\!xtwo  \advance \!dimenA -\!xone
  \!dimenB=\!ytwo  \advance \!dimenB -\!yone
  \ifdim \!dimenA<\!dimenB
    \!startvshade (\!xone,\!yone,\!ytwo)
    \!lshade      (\!xtwo,\!yone,\!ytwo)
  \else
    \!starthshade (\!yone,\!xone,\!xtwo)
    \!lshade      (\!ytwo,\!xone,\!xtwo)
  \fi
  \ignorespaces}
  
\def\frame{%
  \!ifnextchar<{\!frame}{\!frame<\!zpt> }}
\long\def\!frame<#1> #2{%
  \beginpicture
    \setcoordinatesystem units <1pt,1pt> point at 0 0 
    \put {#2} [Bl] at 0 0 
    \!dimenA=#1\relax
    \!dimenB=\!wd \advance \!dimenB \!dimenA
    \!dimenC=\!ht \advance \!dimenC \!dimenA
    \!dimenD=\!dp \advance \!dimenD \!dimenA
    \let\!MFr=\!M
    \!setdimenmode
    \putrectangle corners at {-\!dimenA} {-\!dimenD} and {\!dimenB} {\!dimenC}
    \!setcoordmode
    \let\!M=\!MFr
  \endpicture
  \ignorespaces}
 
\def\rectangle <#1> <#2> {%
  \setbox0=\hbox{}\wd0=#1\ht0=#2\frame {\box0}}

%


\def\!plotfromfile"#1"{%
  \expandafter\!drawcurve \input #1 /}

\def\setquadratic{%
  \let\!drawcurve=\!qcurve
  \let\!!Shade=\!!qShade
  \let\!!!Shade=\!!!qShade}

\def\setlinear{%
  \let\!drawcurve=\!lcurve
  \let\!!Shade=\!!lShade
  \let\!!!Shade=\!!!lShade}

\def\sethistograms{%
  \let\!drawcurve=\!hcurve}

\def\!qcurve #1 #2 {%
  \!start (#1,#2)
  \!Qjoin}
\def\!Qjoin#1 #2 #3 #4 {%
  \!qjoin (#1,#2) (#3,#4)             
  \!ifnextchar/{\!finish}{\!Qjoin}}

\def\!lcurve #1 #2 {%
  \!start (#1,#2)
  \!Ljoin}
\def\!Ljoin#1 #2 {%
  \!ljoin (#1,#2)                    
  \!ifnextchar/{\!finish}{\!Ljoin}}

\def\!finish/{\ignorespaces}

\def\!hcurve #1 #2 {%
  \edef\!hxS{#1}%
  \edef\!hyS{#2}%
  \!hjoin}
\def\!hjoin#1 #2 {%
  \putrectangle corners at {\!hxS} {\!hyS} and {#1} {#2}
  \edef\!hxS{#1}%
  \!ifnextchar/{\!finish}{\!hjoin}}

\def\vshade #1 #2 #3 {%
  \!startvshade (#1,#2,#3)
  \!Shadewhat}

\def\hshade #1 #2 #3 {%
  \!starthshade (#1,#2,#3)
  \!Shadewhat}

\def\!Shadewhat{%
  \futurelet\!nextchar\!Shade}
\def\!Shade{%
  \if <\!nextchar
    \def\!nextShade{\!!Shade}%
  \else
    \if /\!nextchar
      \def\!nextShade{\!finish}%
    \else
      \def\!nextShade{\!!!Shade}%
    \fi
  \fi
  \!nextShade}
\def\!!lShade<#1> #2 #3 #4 {%
  \!lshade <#1> (#2,#3,#4)                 
  \!Shadewhat}
\def\!!!lShade#1 #2 #3 {%
  \!lshade (#1,#2,#3)
  \!Shadewhat} 
\def\!!qShade<#1> #2 #3 #4 #5 #6 #7 {%
  \!qshade <#1> (#2,#3,#4) (#5,#6,#7)      
  \!Shadewhat}
\def\!!!qShade#1 #2 #3 #4 #5 #6 {%
  \!qshade (#1,#2,#3) (#4,#5,#6)
  \!Shadewhat} 

\setlinear

\def\setdashpattern <#1>{%
  \def\!Flist{}\def\!Blist{}\def\!UDlist{}%
  \!countA=0
  \!ecfor\!item:=#1\do{%
    \!dimenA=\!item\relax
    \expandafter\!rightappend\the\!dimenA\withCS{\\}\to\!UDlist%
    \advance\!countA  1
    \ifodd\!countA
      \expandafter\!rightappend\the\!dimenA\withCS{\!Rule}\to\!Flist%
      \expandafter\!leftappend\the\!dimenA\withCS{\!Rule}\to\!Blist%
    \else 
      \expandafter\!rightappend\the\!dimenA\withCS{\!Skip}\to\!Flist%
      \expandafter\!leftappend\the\!dimenA\withCS{\!Skip}\to\!Blist%
    \fi}%
  \!leaderlength=\!zpt
  \def\!Rule##1{\advance\!leaderlength  ##1}%
  \def\!Skip##1{\advance\!leaderlength  ##1}%
  \!Flist%
  \ifdim\!leaderlength>\!zpt 
  \else
    \def\!Flist{\!Skip{24in}}\def\!Blist{\!Skip{24in}}\ignorespaces
    \def\!UDlist{\\{\!zpt}\\{24in}}\ignorespaces
    \!leaderlength=24in
  \fi
  \!dashingon}

\def\!dashingon{%
  \def\!advancedashing{\!!advancedashing}%
  \def\!drawlinearsegment{\!lineardashed}%
  \def\!puthline{\!putdashedhline}%
  \def\!putvline{\!putdashedvline}%
  \ignorespaces}%
\def\!dashingoff{%
  \def\!advancedashing{\relax}%
  \def\!drawlinearsegment{\!linearsolid}%
  \def\!puthline{\!putsolidhline}%
  \def\!putvline{\!putsolidvline}%
  \ignorespaces}

\def\setdots{%
  \!ifnextchar<{\!setdots}{\!setdots<5pt>}}
\def\!setdots<#1>{%
  \!dimenB=#1\advance\!dimenB -\plotsymbolspacing
  \ifdim\!dimenB<\!zpt
    \!dimenB=\!zpt
  \fi
\setdashpattern <\plotsymbolspacing,\!dimenB>}
 
\def\setdotsnear <#1> for <#2>{%
  \!dimenB=#2\relax  \advance\!dimenB -.05pt  
  \!dimenC=#1\relax  \!countA=\!dimenC 
  \!dimenD=\!dimenB  \advance\!dimenD .5\!dimenC  \!countB=\!dimenD
  \divide \!countB  \!countA
  \ifnum 1>\!countB 
    \!countB=1
  \fi
  \divide\!dimenB  \!countB
  \setdots <\!dimenB>}
 
\def\setdashes{%
  \!ifnextchar<{\!setdashes}{\!setdashes<5pt>}}
\def\!setdashes<#1>{\setdashpattern <#1,#1>}
 
\def\setdashesnear <#1> for <#2>{%
  \!dimenB=#2\relax  
  \!dimenC=#1\relax  \!countA=\!dimenC 
  \!dimenD=\!dimenB  \advance\!dimenD .5\!dimenC  \!countB=\!dimenD
  \divide \!countB  \!countA
  \ifodd \!countB 
  \else 
    \advance \!countB  1
  \fi
  \divide\!dimenB  \!countB
  \setdashes <\!dimenB>}
 
\def\setsolid{%
  \def\!Flist{\!Rule{24in}}\def\!Blist{\!Rule{24in}}%
  \def\!UDlist{\\{24in}\\{\!zpt}}%
  \!dashingoff}  
\setsolid


 
  
 
\def\!divide#1#2#3{%
  \!dimenB=#1
  \!dimenC=#2
  \!dimenD=\!dimenB
  \divide \!dimenD \!dimenC
  \!dimenA=\!dimenD
  \multiply\!dimenD \!dimenC
  \advance\!dimenB -\!dimenD
  \!dimenD=\!dimenC
    \ifdim\!dimenD<\!zpt \!dimenD=-\!dimenD 
  \fi
  \ifdim\!dimenD<64pt
    \!divstep[\!tfs]\!divstep[\!tfs]%
  \else 
    \!!divide
  \fi
  #3=\!dimenA\ignorespaces}

\def\!!divide{%
  \ifdim\!dimenD<256pt
    \!divstep[64]\!divstep[32]\!divstep[32]%
  \else 
    \!divstep[8]\!divstep[8]\!divstep[8]\!divstep[8]\!divstep[8]%
    \!dimenA=2\!dimenA
  \fi}

\def\!divstep[#1]{
  \!dimenB=#1\!dimenB
  \!dimenD=\!dimenB
    \divide \!dimenD by \!dimenC
  \!dimenA=#1\!dimenA
    \advance\!dimenA by \!dimenD%
  \multiply\!dimenD by \!dimenC
    \advance\!dimenB by -\!dimenD}
 
\def\Divide <#1> by <#2> forming <#3> {%
  \!divide{#1}{#2}{#3}}

 
 

 

\def\ellipticalarc axes ratio #1:#2 #3 degrees from #4 #5 center at #6 #7 {%
  \!angle=#3pt\relax
  \ifdim\!angle>\!zpt 
    \def\!sign{}
  \else 
    \def\!sign{-}\!angle=-\!angle
  \fi
  \!xxloc=\!M{#6}\!xunit
  \!yyloc=\!M{#7}\!yunit     
  \!xxS=\!M{#4}\!xunit
  \!yyS=\!M{#5}\!yunit
  \advance\!xxS -\!xxloc
  \advance\!yyS -\!yyloc
  \!divide\!xxS{#1pt}\!xxS 
  \!divide\!yyS{#2pt}\!yyS 
  \let\!MC=\!M
  \!setdimenmode
  \!xS=#1\!xxS  \advance\!xS\!xxloc
  \!yS=#2\!yyS  \advance\!yS\!yyloc
  \!start (\!xS,\!yS)%
  \!loop\ifdim\!angle>14.9999pt
    \!rotate(\!xxS,\!yyS)by(\!cos,\!sign\!sin)to(\!xxM,\!yyM) 
    \!rotate(\!xxM,\!yyM)by(\!cos,\!sign\!sin)to(\!xxE,\!yyE)
    \!xM=#1\!xxM  \advance\!xM\!xxloc  \!yM=#2\!yyM  \advance\!yM\!yyloc
    \!xE=#1\!xxE  \advance\!xE\!xxloc  \!yE=#2\!yyE  \advance\!yE\!yyloc
    \!qjoin (\!xM,\!yM) (\!xE,\!yE)
    \!xxS=\!xxE  \!yyS=\!yyE 
    \advance \!angle -15pt
  \repeat
  \ifdim\!angle>\!zpt
    \!angle=100.53096\!angle
    \divide \!angle 360 
    \!sinandcos\!angle\!!sin\!!cos
    \!rotate(\!xxS,\!yyS)by(\!!cos,\!sign\!!sin)to(\!xxM,\!yyM) 
    \!rotate(\!xxM,\!yyM)by(\!!cos,\!sign\!!sin)to(\!xxE,\!yyE)
    \!xM=#1\!xxM  \advance\!xM\!xxloc  \!yM=#2\!yyM  \advance\!yM\!yyloc
    \!xE=#1\!xxE  \advance\!xE\!xxloc  \!yE=#2\!yyE  \advance\!yE\!yyloc
    \!qjoin (\!xM,\!yM) (\!xE,\!yE)
  \fi
  \let\!M=\!MC
  \ignorespaces}

\def\!rotate(#1,#2)by(#3,#4)to(#5,#6){%
  \!dimenA=#3#1\advance \!dimenA -#4#2
  \!dimenB=#3#2\advance \!dimenB  #4#1
  \divide \!dimenA 32  \divide \!dimenB 32 
  #5=\!dimenA  #6=\!dimenB
  \ignorespaces}
\def\!sin{4.17684}
\def\!cos{31.72624}

\def\!sinandcos#1#2#3{%
 \!dimenD=#1
 \!dimenA=\!dimenD
 \!dimenB=32pt
 \!removept\!dimenD\!value
 \!dimenC=\!dimenD
 \!dimenC=\!value\!dimenC \divide\!dimenC by 64 
 \advance\!dimenB by -\!dimenC
 \!dimenC=\!value\!dimenC \divide\!dimenC by 96 
 \advance\!dimenA by -\!dimenC
 \!dimenC=\!value\!dimenC \divide\!dimenC by 128 
 \advance\!dimenB by \!dimenC%
 \!removept\!dimenA#2
 \!removept\!dimenB#3
 \ignorespaces}




\def\putrule#1from #2 #3 to #4 #5 {%
  \!xloc=\!M{#2}\!xunit  \!xxloc=\!M{#4}\!xunit%
  \!yloc=\!M{#3}\!yunit  \!yyloc=\!M{#5}\!yunit%
  \!dxpos=\!xxloc  \advance\!dxpos by -\!xloc
  \!dypos=\!yyloc  \advance\!dypos by -\!yloc
  \ifdim\!dypos=\!zpt
    \def\!!Line{\!puthline{#1}}\ignorespaces
  \else
    \ifdim\!dxpos=\!zpt
      \def\!!Line{\!putvline{#1}}\ignorespaces
    \else 
       \def\!!Line{}
    \fi
  \fi
  \let\!ML=\!M
  \!setdimenmode
  \!!Line%
  \let\!M=\!ML
  \ignorespaces}

\def\!putsolidhline#1{%
  \ifdim\!dxpos>\!zpt 
    \put{\!hline\!dxpos}#1[l] at {\!xloc} {\!yloc}
  \else 
    \put{\!hline{-\!dxpos}}#1[l] at {\!xxloc} {\!yyloc}
  \fi
  \ignorespaces}
 
\def\!putsolidvline#1{%
  \ifdim\!dypos>\!zpt 
    \put{\!vline\!dypos}#1[b] at {\!xloc} {\!yloc}
  \else 
    \put{\!vline{-\!dypos}}#1[b] at {\!xxloc} {\!yyloc}
  \fi
  \ignorespaces}
 
\def\!hline#1{\hbox to #1{\leaders \hrule height\linethickness\hfill}}
\def\!vline#1{\vbox to #1{\leaders \vrule width\linethickness\vfill}}

\def\!putdashedhline#1{%
  \ifdim\!dxpos>\!zpt 
    \!DLsetup\!Flist\!dxpos
    \put{\hbox to \!totalleaderlength{\!hleaders}\!hpartialpattern\!Rtrunc}
      #1[l] at {\!xloc} {\!yloc} 
  \else 
    \!DLsetup\!Blist{-\!dxpos}
    \put{\!hpartialpattern\!Ltrunc\hbox to \!totalleaderlength{\!hleaders}}
      #1[r] at {\!xloc} {\!yloc} 
  \fi
  \ignorespaces}
 
\def\!putdashedvline#1{%
  \!dypos=-\!dypos
  \ifdim\!dypos>\!zpt 
    \!DLsetup\!Flist\!dypos 
    \put{\vbox{\vbox to \!totalleaderlength{\!vleaders}
      \!vpartialpattern\!Rtrunc}}#1[t] at {\!xloc} {\!yloc} 
  \else 
    \!DLsetup\!Blist{-\!dypos}
    \put{\vbox{\!vpartialpattern\!Ltrunc
      \vbox to \!totalleaderlength{\!vleaders}}}#1[b] at {\!xloc} {\!yloc} 
  \fi
  \ignorespaces}

\def\!DLsetup#1#2{
  \let\!RSlist=#1
  \!countB=#2
  \!countA=\!leaderlength
  \divide\!countB by \!countA
  \!totalleaderlength=\!countB\!leaderlength
  \!Rresiduallength=#2%
  \advance \!Rresiduallength by -\!totalleaderlength
  \!Lresiduallength=\!leaderlength
  \advance \!Lresiduallength by -\!Rresiduallength
  \ignorespaces}
 
\def\!hleaders{%
  \def\!Rule##1{\vrule height\linethickness width##1}%
  \def\!Skip##1{\hskip##1}%
  \leaders\hbox{\!RSlist}\hfill}
 
\def\!hpartialpattern#1{%
  \!dimenA=\!zpt \!dimenB=\!zpt 
  \def\!Rule##1{#1{##1}\vrule height\linethickness width\!dimenD}%
  \def\!Skip##1{#1{##1}\hskip\!dimenD}%
  \!RSlist}
 
\def\!vleaders{%
  \def\!Rule##1{\hrule width\linethickness height##1}%
  \def\!Skip##1{\vskip##1}%
  \leaders\vbox{\!RSlist}\vfill}
 
\def\!vpartialpattern#1{%
  \!dimenA=\!zpt \!dimenB=\!zpt 
  \def\!Rule##1{#1{##1}\hrule width\linethickness height\!dimenD}%
  \def\!Skip##1{#1{##1}\vskip\!dimenD}%
  \!RSlist}
 
\def\!Rtrunc#1{\!trunc{#1}>\!Rresiduallength}
\def\!Ltrunc#1{\!trunc{#1}<\!Lresiduallength}
 
\def\!trunc#1#2#3{%
  \!dimenA=\!dimenB         
  \advance\!dimenB by #1%
  \!dimenD=\!dimenB  \ifdim\!dimenD#2#3\!dimenD=#3\fi
  \!dimenC=\!dimenA  \ifdim\!dimenC#2#3\!dimenC=#3\fi
  \advance \!dimenD by -\!dimenC}

\def\!start (#1,#2){%
  \!plotxorigin=\!xorigin  \advance \!plotxorigin by \!plotsymbolxshift
  \!plotyorigin=\!yorigin  \advance \!plotyorigin by \!plotsymbolyshift
  \!xS=\!M{#1}\!xunit \!yS=\!M{#2}\!yunit
  \!rotateaboutpivot\!xS\!yS
  \!copylist\!UDlist\to\!!UDlist
  \!getnextvalueof\!downlength\from\!!UDlist
  \!distacross=\!zpt
  \!intervalno=0 
  \global\totalarclength=\!zpt
  \ignorespaces}

\def\!ljoin (#1,#2){%
  \advance\!intervalno by 1
  \!xE=\!M{#1}\!xunit \!yE=\!M{#2}\!yunit
  \!rotateaboutpivot\!xE\!yE
  \!xdiff=\!xE \advance \!xdiff by -\!xS
  \!ydiff=\!yE \advance \!ydiff by -\!yS
  \!Pythag\!xdiff\!ydiff\!arclength
  \global\advance \totalarclength by \!arclength%
  \!drawlinearsegment
  \!xS=\!xE \!yS=\!yE
  \ignorespaces}

\def\!linearsolid{%
  \!npoints=\!arclength
  \!countA=\plotsymbolspacing
  \divide\!npoints by \!countA
  \ifnum \!npoints<1 
    \!npoints=1 
  \fi
  \divide\!xdiff by \!npoints
  \divide\!ydiff by \!npoints
  \!xpos=\!xS \!ypos=\!yS
  \loop\ifnum\!npoints>-1
    \!plotifinbounds
    \advance \!xpos by \!xdiff
    \advance \!ypos by \!ydiff
    \advance \!npoints by -1
  \repeat
  \ignorespaces}

\def\!lineardashed{%
  \ifdim\!distacross>\!arclength
    \advance \!distacross by -\!arclength  
  \else
    \loop\ifdim\!distacross<\!arclength
      \!divide\!distacross\!arclength\!dimenA
      \!removept\!dimenA\!t
      \!xpos=\!t\!xdiff \advance \!xpos by \!xS
      \!ypos=\!t\!ydiff \advance \!ypos by \!yS
      \!plotifinbounds
      \advance\!distacross by \plotsymbolspacing
      \!advancedashing
    \repeat  
    \advance \!distacross by -\!arclength
  \fi
  \ignorespaces}

\def\!!advancedashing{%
  \advance\!downlength by -\plotsymbolspacing
  \ifdim \!downlength>\!zpt
  \else
    \advance\!distacross by \!downlength
    \!getnextvalueof\!uplength\from\!!UDlist
    \advance\!distacross by \!uplength
    \!getnextvalueof\!downlength\from\!!UDlist
  \fi}

\def\inboundscheckoff{%
  \def\!plotifinbounds{\!plot(\!xpos,\!ypos)}%
  \def\!initinboundscheck{\relax}\ignorespaces}
 
\inboundscheckoff
 
\def\!!plotifinbounds{%
  \ifdim \!xpos<\!checkleft
  \else
    \ifdim \!xpos>\!checkright
    \else
      \ifdim \!ypos<\!checkbot
      \else
         \ifdim \!ypos>\!checktop
         \else
           \!plot(\!xpos,\!ypos)
         \fi 
      \fi
    \fi
  \fi}

\def\!!initinboundscheck{%
  \!checkleft=\!arealloc     \advance\!checkleft by \!xorigin
  \!checkright=\!arearloc    \advance\!checkright by \!xorigin
  \!checkbot=\!areabloc      \advance\!checkbot by \!yorigin
  \!checktop=\!areatloc      \advance\!checktop by \!yorigin}

%


\def\!logten#1#2{%
  \expandafter\!!logten#1\!nil
  \!removept\!dimenF#2%
  \ignorespaces}

\def\!!logten#1#2\!nil{%
  \if -#1%
    \!dimenF=\!zpt
    \def\!next{\ignorespaces}%
  \else
    \if +#1%
      \def\!next{\!!logten#2\!nil}%
    \else
      \if .#1%
        \def\!next{\!!logten0.#2\!nil}%
      \else
        \def\!next{\!!!logten#1#2..\!nil}%
      \fi
    \fi
  \fi
  \!next}

\def\!!!logten#1#2.#3.#4\!nil{%
  \!dimenF=1pt 
  \if 0#1%
    \!!logshift#3pt 
  \else 
    \!logshift#2/
    \!dimenE=#1.#2#3pt 
  \fi 
  \ifdim \!dimenE<\!rootten
    \multiply \!dimenE 10 
    \advance  \!dimenF -1pt
  \fi
  \!dimenG=\!dimenE
    \advance\!dimenG 10pt
  \advance\!dimenE -10pt 
  \multiply\!dimenE 10 
  \!divide\!dimenE\!dimenG\!dimenE
  \!removept\!dimenE\!t
  \!dimenG=\!t\!dimenE
  \!removept\!dimenG\!tt
  \!dimenH=\!tt\!tenAe
    \divide\!dimenH 100
  \advance\!dimenH \!tenAc
  \!dimenH=\!tt\!dimenH
    \divide\!dimenH 100   
  \advance\!dimenH \!tenAa
  \!dimenH=\!t\!dimenH
    \divide\!dimenH 100 
  \advance\!dimenF \!dimenH}

\def\!logshift#1{%
  \if #1/%
    \def\!next{\ignorespaces}%
  \else
    \advance\!dimenF 1pt 
    \def\!next{\!logshift}%
  \fi 
  \!next}
 
 \def\!!logshift#1{%
   \advance\!dimenF -1pt
   \if 0#1%
     \def\!next{\!!logshift}%
   \else
     \if p#1%
       \!dimenF=1pt
       \def\!next{\!dimenE=1p}%
     \else
       \def\!next{\!dimenE=#1.}%
     \fi
   \fi
   \!next}

\def\beginpicture{%
  \setbox\!picbox=\hbox\bgroup%
  \!xleft=\maxdimen  
  \!xright=-\maxdimen
  \!ybot=\maxdimen
  \!ytop=-\maxdimen}
 
\def\endpicture{%
  \ifdim\!xleft=\maxdimen
    \!xleft=\!zpt \!xright=\!zpt \!ybot=\!zpt \!ytop=\!zpt 
  \fi
  \global\!Xleft=\!xleft \global\!Xright=\!xright
  \global\!Ybot=\!ybot \global\!Ytop=\!ytop
  \egroup%
  \ht\!picbox=\!Ytop  \dp\!picbox=-\!Ybot
  \ifdim\!Ybot>\!zpt
  \else 
    \ifdim\!Ytop<\!zpt
      \!Ybot=\!Ytop
    \else
      \!Ybot=\!zpt
    \fi
  \fi
  \hbox{\kern-\!Xleft\lower\!Ybot\box\!picbox\kern\!Xright}}
 
\def\endpicturesave <#1,#2>{%
  \endpicture \global #1=\!Xleft \global #2=\!Ybot \ignorespaces}

\def\setcoordinatesystem{%
  \!ifnextchar{u}{\!getlengths }
    {\!getlengths units <\!xunit,\!yunit>}}
\def\!getlengths units <#1,#2>{%
  \!xunit=#1\relax
  \!yunit=#2\relax
  \!ifcoordmode 
    \let\!SCnext=\!SCccheckforRP
  \else
    \let\!SCnext=\!SCdcheckforRP
  \fi
  \!SCnext}
\def\!SCccheckforRP{%
  \!ifnextchar{p}{\!cgetreference }
    {\!cgetreference point at {\!xref} {\!yref} }}
\def\!cgetreference point at #1 #2 {%
  \edef\!xref{#1}\edef\!yref{#2}%
  \!xorigin=\!xref\!xunit  \!yorigin=\!yref\!yunit  
  \!initinboundscheck 
  \ignorespaces}
\def\!SCdcheckforRP{%
  \!ifnextchar{p}{\!dgetreference}%
    {\ignorespaces}}
\def\!dgetreference point at #1 #2 {%
  \!xorigin=#1\relax  \!yorigin=#2\relax
  \ignorespaces}

\long\def\put#1#2 at #3 #4 {%
  \!setputobject{#1}{#2}%
  \!xpos=\!M{#3}\!xunit  \!ypos=\!M{#4}\!yunit  
  \!rotateaboutpivot\!xpos\!ypos%
  \advance\!xpos -\!xorigin  \advance\!xpos -\!xshift
  \advance\!ypos -\!yorigin  \advance\!ypos -\!yshift
  \kern\!xpos\raise\!ypos\box\!putobject\kern-\!xpos%
  \!doaccounting\ignorespaces}
 
\long\def\multiput #1#2 at {%
  \!setputobject{#1}{#2}%
  \!ifnextchar"{\!putfromfile}{\!multiput}}
\def\!putfromfile"#1"{%
  \expandafter\!multiput \input #1 /}
\def\!multiput{%
  \futurelet\!nextchar\!!multiput}
\def\!!multiput{%
  \if *\!nextchar
    \def\!nextput{\!alsoby}%
  \else
    \if /\!nextchar
      \def\!nextput{\!finishmultiput}%
    \else
      \def\!nextput{\!alsoat}%
    \fi
  \fi
  \!nextput}
\def\!finishmultiput/{%
  \setbox\!putobject=\hbox{}%
  \ignorespaces}
 
\def\!alsoat#1 #2 {%
  \!xpos=\!M{#1}\!xunit  \!ypos=\!M{#2}\!yunit  
  \!rotateaboutpivot\!xpos\!ypos%
  \advance\!xpos -\!xorigin  \advance\!xpos -\!xshift
  \advance\!ypos -\!yorigin  \advance\!ypos -\!yshift
  \kern\!xpos\raise\!ypos\copy\!putobject\kern-\!xpos%
  \!doaccounting
  \!multiput}
 
\def\!alsoby*#1 #2 #3 {%
  \!dxpos=\!M{#2}\!xunit \!dypos=\!M{#3}\!yunit 
  \!rotateonly\!dxpos\!dypos
  \!ntemp=#1%
  \!!loop\ifnum\!ntemp>0
    \advance\!xpos by \!dxpos  \advance\!ypos by \!dypos
    \kern\!xpos\raise\!ypos\copy\!putobject\kern-\!xpos%
    \advance\!ntemp by -1
  \repeat
  \!doaccounting 
  \!multiput}
 
\def\accountingon{\def\!doaccounting{\!!doaccounting}\ignorespaces}

\accountingon
\def\!!doaccounting{%
  \!xtemp=\!xpos  
  \!ytemp=\!ypos
  \ifdim\!xtemp<\!xleft 
     \!xleft=\!xtemp 
  \fi
  \advance\!xtemp by  \!wd 
  \ifdim\!xright<\!xtemp 
    \!xright=\!xtemp
  \fi
  \advance\!ytemp by -\!dp
  \ifdim\!ytemp<\!ybot  
    \!ybot=\!ytemp
  \fi
  \advance\!ytemp by  \!dp
  \advance\!ytemp by  \!ht 
  \ifdim\!ytemp>\!ytop  
    \!ytop=\!ytemp  
  \fi}
 
\long\def\!setputobject#1#2{%
  \setbox\!putobject=\hbox{#1}%
  \!ht=\ht\!putobject  \!dp=\dp\!putobject  \!wd=\wd\!putobject
  \wd\!putobject=\!zpt
  \!xshift=.5\!wd   \!yshift=.5\!ht   \advance\!yshift by -.5\!dp
  \edef\!putorientation{#2}%
  \expandafter\!SPOreadA\!putorientation[]\!nil%
  \expandafter\!SPOreadB\!putorientation<\!zpt,\!zpt>\!nil\ignorespaces}
 
\def\!SPOreadA#1[#2]#3\!nil{\!etfor\!orientation:=#2\do\!SPOreviseshift}
 
\def\!SPOreadB#1<#2,#3>#4\!nil{\advance\!xshift by -#2\advance\!yshift by -#3}
 
\def\!SPOreviseshift{%
  \if l\!orientation 
    \!xshift=\!zpt
  \else 
    \if r\!orientation 
      \!xshift=\!wd
    \else 
      \if b\!orientation
        \!yshift=-\!dp
      \else 
        \if B\!orientation 
          \!yshift=\!zpt
        \else 
          \if t\!orientation 
            \!yshift=\!ht
          \fi 
        \fi
      \fi
    \fi
  \fi}

\long\def\!dimenput#1#2(#3,#4){%
  \!setputobject{#1}{#2}%
  \!xpos=#3\advance\!xpos by -\!xshift
  \!ypos=#4\advance\!ypos by -\!yshift
  \kern\!xpos\raise\!ypos\box\!putobject\kern-\!xpos%
  \!doaccounting\ignorespaces}

\def\!setdimenmode{%
  \let\!M=\!M!!\ignorespaces}
\def\!setcoordmode{%
  \let\!M=\!M!\ignorespaces}
\def\!ifcoordmode{%
  \ifx \!M \!M!}
\def\!ifdimenmode{%
  \ifx \!M \!M!!}
\def\!M!#1#2{#1#2} 
\def\!M!!#1#2{#1}
\!setcoordmode
\let\setdimensionmode=\!setdimenmode
\let\setcoordinatemode=\!setcoordmode




\def\!stack[#1]{%
  \let\!lglue=\hfill \let\!rglue=\hfill
  \expandafter\let\csname !#1glue\endcsname=\relax
  \!ifnextchar<{\!!stack}{\!!stack<\stackleading>}}
\def\!!stack<#1>#2{%
  \vbox{\def\!valueslist{}\!ecfor\!value:=#2\do{%
    \expandafter\!rightappend\!value\withCS{\\}\to\!valueslist}%
    \!lop\!valueslist\to\!value
    \let\\=\cr\lineskiplimit=\maxdimen\lineskip=#1%
    \baselineskip=-1000pt\halign{\!lglue##\!rglue\cr \!value\!valueslist\cr}}%
  \ignorespaces}


\def\!lines[#1]#2{%
  \let\!lglue=\hfill \let\!rglue=\hfill
  \expandafter\let\csname !#1glue\endcsname=\relax
  \vbox{\halign{\!lglue##\!rglue\cr #2\crcr}}%
  \ignorespaces}


\def\!Lines[#1]#2{%
  \let\!lglue=\hfill \let\!rglue=\hfill
  \expandafter\let\csname !#1glue\endcsname=\relax
  \vtop{\halign{\!lglue##\!rglue\cr #2\crcr}}%
  \ignorespaces}

 
 
 
\def\setplotsymbol(#1#2){%
  \!setputobject{#1}{#2}
  \setbox\!plotsymbol=\box\!putobject%
  \!plotsymbolxshift=\!xshift 
  \!plotsymbolyshift=\!yshift 
  \ignorespaces}
 
\setplotsymbol({\fiverm .})

 
\def\!!plot(#1,#2){%
  \!dimenA=-\!plotxorigin \advance \!dimenA by #1
  \!dimenB=-\!plotyorigin \advance \!dimenB by #2
  \kern\!dimenA\raise\!dimenB\copy\!plotsymbol\kern-\!dimenA%
  \ignorespaces}
 
\def\!!!plot(#1,#2){%
  \!dimenA=-\!plotxorigin \advance \!dimenA by #1
  \!dimenB=-\!plotyorigin \advance \!dimenB by #2
  \kern\!dimenA\raise\!dimenB\copy\!plotsymbol\kern-\!dimenA%
  \!countE=\!dimenA
  \!countF=\!dimenB
  \immediate\write\!replotfile{\the\!countE,\the\!countF.}%
  \ignorespaces}

\def\savelinesandcurves on "#1" {%
  \immediate\closeout\!replotfile
  \immediate\openout\!replotfile=#1%
  \let\!plot=\!!!plot}

\def\dontsavelinesandcurves {%
  \let\!plot=\!!plot}
\dontsavelinesandcurves

{\catcode`\%=11\xdef\!Commentsignal{
\def\writesavefile#1 {%
  \immediate\write\!replotfile{\!Commentsignal #1}%
  \ignorespaces}

\def\replot"#1" {%
  \expandafter\!replot\input #1 /}
\def\!replot#1,#2. {%
  \!dimenA=#1sp
  \kern\!dimenA\raise#2sp\copy\!plotsymbol\kern-\!dimenA
  \futurelet\!nextchar\!!replot}
\def\!!replot{%
  \if /\!nextchar 
    \def\!next{\!finish}%
  \else
    \def\!next{\!replot}%
  \fi
  \!next}


 
 
\def\!Pythag#1#2#3{%
  \!dimenE=#1\relax                                     
  \ifdim\!dimenE<\!zpt 
    \!dimenE=-\!dimenE 
  \fi
  \!dimenF=#2\relax
  \ifdim\!dimenF<\!zpt 
    \!dimenF=-\!dimenF 
  \fi
  \advance \!dimenF by \!dimenE
  \ifdim\!dimenF=\!zpt 
    \!dimenG=\!zpt
  \else 
    \!divide{8\!dimenE}\!dimenF\!dimenE
    \advance\!dimenE by -4pt
      \!dimenE=2\!dimenE
    \!removept\!dimenE\!!t
    \!dimenE=\!!t\!dimenE
    \advance\!dimenE by 64pt
    \divide \!dimenE by 2
    \!dimenH=7pt
    \!!Pythag\!!Pythag\!!Pythag
    \!removept\!dimenH\!!t
    \!dimenG=\!!t\!dimenF
    \divide\!dimenG by 8
  \fi
  #3=\!dimenG
  \ignorespaces}

\def\!!Pythag{
  \!divide\!dimenE\!dimenH\!dimenI
  \advance\!dimenH by \!dimenI
    \divide\!dimenH by 2}

\def\placehypotenuse for <#1> and <#2> in <#3> {%
  \!Pythag{#1}{#2}{#3}}

 
 
 
\def\!qjoin (#1,#2) (#3,#4){%
  \advance\!intervalno by 1
  \!ifcoordmode
    \edef\!xmidpt{#1}\edef\!ymidpt{#2}%
  \else
    \!dimenA=#1\relax \edef\!xmidpt{\the\!dimenA}%
    \!dimenA=#2\relax \edef\!xmidpt{\the\!dimenA}%
  \fi
  \!xM=\!M{#1}\!xunit  \!yM=\!M{#2}\!yunit   \!rotateaboutpivot\!xM\!yM
  \!xE=\!M{#3}\!xunit  \!yE=\!M{#4}\!yunit   \!rotateaboutpivot\!xE\!yE
%
  \!dimenA=\!xM  \advance \!dimenA by -\!xS
  \!dimenB=\!xE  \advance \!dimenB by -\!xM
  \!xB=3\!dimenA \advance \!xB by -\!dimenB
  \!xC=2\!dimenB \advance \!xC by -2\!dimenA
%
  \!dimenA=\!yM  \advance \!dimenA by -\!yS%
  \!dimenB=\!yE  \advance \!dimenB by -\!yM%
  \!yB=3\!dimenA \advance \!yB by -\!dimenB%
  \!yC=2\!dimenB \advance \!yC by -2\!dimenA%
%
  \!xprime=\!xB  \!yprime=\!yB
  \!dxprime=.5\!xC  \!dyprime=.5\!yC
  \!getf \!midarclength=\!dimenA
  \!getf \advance \!midarclength by 4\!dimenA
  \!getf \advance \!midarclength by \!dimenA
  \divide \!midarclength by 12
%
  \!arclength=\!dimenA
  \!getf \advance \!arclength by 4\!dimenA
  \!getf \advance \!arclength by \!dimenA
  \divide \!arclength by 12
  \advance \!arclength by \!midarclength
  \global\advance \totalarclength by \!arclength
%
%
  \ifdim\!distacross>\!arclength 
    \advance \!distacross by -\!arclength
  \else
    \!initinverseinterp
    \loop\ifdim\!distacross<\!arclength
      \!inverseinterp
      \!xpos=\!t\!xC \advance\!xpos by \!xB
        \!xpos=\!t\!xpos \advance \!xpos by \!xS
      \!ypos=\!t\!yC \advance\!ypos by \!yB
        \!ypos=\!t\!ypos \advance \!ypos by \!yS
      \!plotifinbounds
      \advance\!distacross \plotsymbolspacing
      \!advancedashing
    \repeat  
    \advance \!distacross by -\!arclength
  \fi
  \!xS=\!xE
  \!yS=\!yE
  \ignorespaces}

\def\!getf{\!Pythag\!xprime\!yprime\!dimenA%
  \advance\!xprime by \!dxprime
  \advance\!yprime by \!dyprime}

\def\!initinverseinterp{%
  \ifdim\!arclength>\!zpt
    \!divide{8\!midarclength}\!arclength\!dimenE
    \ifdim\!dimenE<\!wmin \!setinverselinear
    \else 
      \ifdim\!dimenE>\!wmax \!setinverselinear
      \else
        \def\!inverseinterp{\!inversequad}\ignorespaces
%
%
         \!removept\!dimenE\!Ew
         \!dimenF=-\!Ew\!dimenE
         \advance\!dimenF by 32pt
         \!dimenG=8pt 
         \advance\!dimenG by -\!dimenE
         \!dimenG=\!Ew\!dimenG
         \!divide\!dimenF\!dimenG\!beta
         \!gamma=1pt
         \advance \!gamma by -\!beta
      \fi
    \fi
  \fi
  \ignorespaces}

\def\!inversequad{%
  \!divide\!distacross\!arclength\!dimenG
  \!removept\!dimenG\!v
  \!dimenG=\!v\!gamma
  \advance\!dimenG by \!beta
  \!dimenG=\!v\!dimenG
  \!removept\!dimenG\!t}

\def\!setinverselinear{%
  \def\!inverseinterp{\!inverselinear}%
  \divide\!dimenE by 8 \!removept\!dimenE\!t
  \!countC=\!intervalno \multiply \!countC 2
  \!countB=\!countC     \advance \!countB -1
  \!countA=\!countB     \advance \!countA -1
  \wlog{\the\!countB th point (\!xmidpt,\!ymidpt) being plotted 
    doesn't lie in the}%
  \wlog{ middle third of the arc between the \the\!countA th 
    and \the\!countC th points:}%
  \wlog{ [arc length \the\!countA\space to \the\!countB]/[arc length 
    \the \!countA\space to \the\!countC]=\!t.}%
  \ignorespaces}
 
\def\!inverselinear{%
  \!divide\!distacross\!arclength\!dimenG
  \!removept\!dimenG\!t}

 

\def\startrotation{%
  \let\!rotateaboutpivot=\!!rotateaboutpivot
  \let\!rotateonly=\!!rotateonly
  \!ifnextchar{b}{\!getsincos }%
    {\!getsincos by {\!cosrotationangle} {\!sinrotationangle} }}
\def\!getsincos by #1 #2 {%
  \edef\!cosrotationangle{#1}%
  \edef\!sinrotationangle{#2}%
  \!ifcoordmode 
    \let\!ROnext=\!ccheckforpivot
  \else
    \let\!ROnext=\!dcheckforpivot
  \fi
  \!ROnext}
\def\!ccheckforpivot{%
  \!ifnextchar{a}{\!cgetpivot}%
    {\!cgetpivot about {\!xpivotcoord} {\!ypivotcoord} }}
\def\!cgetpivot about #1 #2 {%
  \edef\!xpivotcoord{#1}%
  \edef\!ypivotcoord{#2}%
  \!xpivot=#1\!xunit  \!ypivot=#2\!yunit
  \ignorespaces}
\def\!dcheckforpivot{%
  \!ifnextchar{a}{\!dgetpivot}{\ignorespaces}}
\def\!dgetpivot about #1 #2 {%
  \!xpivot=#1\relax  \!ypivot=#2\relax
  \ignorespaces}

\def\stoprotation{%
  \let\!rotateaboutpivot=\!!!rotateaboutpivot
  \let\!rotateonly=\!!!rotateonly
  \ignorespaces}
 
\def\!!rotateaboutpivot#1#2{%
  \!dimenA=#1\relax  \advance\!dimenA -\!xpivot
  \!dimenB=#2\relax  \advance\!dimenB -\!ypivot
  \!dimenC=\!cosrotationangle\!dimenA
    \advance \!dimenC -\!sinrotationangle\!dimenB
  \!dimenD=\!cosrotationangle\!dimenB
    \advance \!dimenD  \!sinrotationangle\!dimenA
  \advance\!dimenC \!xpivot  \advance\!dimenD \!ypivot
  #1=\!dimenC  #2=\!dimenD
  \ignorespaces}

\def\!!rotateonly#1#2{%
  \!dimenA=#1\relax  \!dimenB=#2\relax 
  \!dimenC=\!cosrotationangle\!dimenA
    \advance \!dimenC -\!rotsign\!sinrotationangle\!dimenB
  \!dimenD=\!cosrotationangle\!dimenB
    \advance \!dimenD  \!rotsign\!sinrotationangle\!dimenA
  #1=\!dimenC  #2=\!dimenD
  \ignorespaces}
\def\!rotsign{}
\def\!!!rotateaboutpivot#1#2{\relax}
\def\!!!rotateonly#1#2{\relax}
\stoprotation

\def\!reverserotateonly#1#2{%
  \def\!rotsign{-}%
  \!rotateonly{#1}{#2}%
  \def\!rotsign{}%
  \ignorespaces}

\def\!getspan span <#1>{%
  \!dshade=#1\relax
  \!ifcoordmode 
    \let\!GRnext=\!GRccheckforAP
  \else
    \let\!GRnext=\!GRdcheckforAP
  \fi
  \!GRnext}
\def\!GRccheckforAP{%
  \!ifnextchar{p}{\!cgetanchor }
    {\!cgetanchor point at {\!xshadesave} {\!yshadesave} }}
\def\!cgetanchor point at #1 #2 {%
  \edef\!xshadesave{#1}\edef\!yshadesave{#2}%
  \!xshade=\!xshadesave\!xunit  \!yshade=\!yshadesave\!yunit
  \ignorespaces}
\def\!GRdcheckforAP{%
  \!ifnextchar{p}{\!dgetanchor}%
    {\ignorespaces}}
\def\!dgetanchor point at #1 #2 {%
  \!xshade=#1\relax  \!yshade=#2\relax
  \ignorespaces}

\def\setshadesymbol{%
  \!ifnextchar<{\!setshadesymbol}{\!setshadesymbol<,,,> }}

\def\!setshadesymbol <#1,#2,#3,#4> (#5#6){%
  \!setputobject{#5}{#6}%
  \setbox\!shadesymbol=\box\!putobject%
  \!shadesymbolxshift=\!xshift \!shadesymbolyshift=\!yshift
%
  \!dimenA=\!xshift \advance\!dimenA \!smidge
  \!override\!dimenA{#1}\!lshrinkage%
  \!dimenA=\!wd \advance \!dimenA -\!xshift
    \advance\!dimenA \!smidge
    \!override\!dimenA{#2}\!rshrinkage
  \!dimenA=\!dp \advance \!dimenA \!yshift
    \advance\!dimenA \!smidge
    \!override\!dimenA{#3}\!bshrinkage
  \!dimenA=\!ht \advance \!dimenA -\!yshift
    \advance\!dimenA \!smidge
    \!override\!dimenA{#4}\!tshrinkage
  \ignorespaces}
\def\!smidge{-.2pt}%

\def\!override#1#2#3{%
  \edef\!!override{#2}%
  \ifx \!!override\empty
    #3=#1\relax
  \else
    \if z\!!override
      #3=\!zpt
    \else
      \ifx \!!override\!blankz
        #3=\!zpt
      \else
        #3=#2\relax
      \fi
    \fi
  \fi
  \ignorespaces}
\def\!blankz{ z}

\setshadesymbol ({\fiverm .})

\def\!startvshade#1(#2,#3,#4){%
  \let\!!xunit=\!xunit%
  \let\!!yunit=\!yunit%
  \let\!!xshade=\!xshade%
  \let\!!yshade=\!yshade%
  \def\!getshrinkages{\!vgetshrinkages}%
  \let\!setshadelocation=\!vsetshadelocation%
  \!xS=\!M{#2}\!!xunit
  \!ybS=\!M{#3}\!!yunit
  \!ytS=\!M{#4}\!!yunit
  \!shadexorigin=\!xorigin  \advance \!shadexorigin \!shadesymbolxshift
  \!shadeyorigin=\!yorigin  \advance \!shadeyorigin \!shadesymbolyshift
  \ignorespaces}
 
\def\!starthshade#1(#2,#3,#4){%
  \let\!!xunit=\!yunit%
  \let\!!yunit=\!xunit%
  \let\!!xshade=\!yshade%
  \let\!!yshade=\!xshade%
  \def\!getshrinkages{\!hgetshrinkages}%
  \let\!setshadelocation=\!hsetshadelocation%
  \!xS=\!M{#2}\!!xunit
  \!ybS=\!M{#3}\!!yunit
  \!ytS=\!M{#4}\!!yunit
  \!shadexorigin=\!xorigin  \advance \!shadexorigin \!shadesymbolxshift
  \!shadeyorigin=\!yorigin  \advance \!shadeyorigin \!shadesymbolyshift
  \ignorespaces}

\def\!lattice#1#2#3#4#5{%
  \!dimenA=#1
  \!dimenB=#2
  \!countB=\!dimenB
%
  \!dimenC=#3
  \advance\!dimenC -\!dimenA
  \!countA=\!dimenC
  \divide\!countA \!countB
  \ifdim\!dimenC>\!zpt
    \!dimenD=\!countA\!dimenB
    \ifdim\!dimenD<\!dimenC
      \advance\!countA 1 
    \fi
  \fi
  \!dimenC=\!countA\!dimenB
    \advance\!dimenC \!dimenA
  #4=\!countA
  #5=\!dimenC
  \ignorespaces}

\def\!qshade#1(#2,#3,#4)#5(#6,#7,#8){%
  \!xM=\!M{#2}\!!xunit
  \!ybM=\!M{#3}\!!yunit
  \!ytM=\!M{#4}\!!yunit
  \!xE=\!M{#6}\!!xunit
  \!ybE=\!M{#7}\!!yunit
  \!ytE=\!M{#8}\!!yunit
  \!getcoeffs\!xS\!ybS\!xM\!ybM\!xE\!ybE\!ybB\!ybC
  \!getcoeffs\!xS\!ytS\!xM\!ytM\!xE\!ytE\!ytB\!ytC
  \def\!getylimits{\!qgetylimits}%
  \!shade{#1}\ignorespaces}
 
\def\!lshade#1(#2,#3,#4){%
  \!xE=\!M{#2}\!!xunit
  \!ybE=\!M{#3}\!!yunit
  \!ytE=\!M{#4}\!!yunit
  \!dimenE=\!xE  \advance \!dimenE -\!xS
  \!dimenC=\!ytE \advance \!dimenC -\!ytS
  \!divide\!dimenC\!dimenE\!ytB
  \!dimenC=\!ybE \advance \!dimenC -\!ybS
  \!divide\!dimenC\!dimenE\!ybB
  \def\!getylimits{\!lgetylimits}%
  \!shade{#1}\ignorespaces}
 
\def\!getcoeffs#1#2#3#4#5#6#7#8{%
  \!dimenC=#4\advance \!dimenC -#2
  \!dimenE=#3\advance \!dimenE -#1
  \!divide\!dimenC\!dimenE\!dimenF
  \!dimenC=#6\advance \!dimenC -#4
  \!dimenH=#5\advance \!dimenH -#3
  \!divide\!dimenC\!dimenH\!dimenG
  \advance\!dimenG -\!dimenF
  \advance \!dimenH \!dimenE
  \!divide\!dimenG\!dimenH#8
  \!removept#8\!t
  #7=-\!t\!dimenE
  \advance #7\!dimenF
  \ignorespaces}

\def\!shade#1{%
  \!getshrinkages#1<,,,>\!nil
  \advance \!dimenE \!xS
  \!lattice\!!xshade\!dshade\!dimenE
    \!parity\!xpos
  \!dimenF=-\!dimenF
    \advance\!dimenF \!xE
  \!loop\!not{\ifdim\!xpos>\!dimenF}
    \!shadecolumn%
    \advance\!xpos \!dshade
    \advance\!parity 1
  \repeat
  \!xS=\!xE
  \!ybS=\!ybE
  \!ytS=\!ytE
  \ignorespaces}

\def\!vgetshrinkages#1<#2,#3,#4,#5>#6\!nil{%
  \!override\!lshrinkage{#2}\!dimenE
  \!override\!rshrinkage{#3}\!dimenF
  \!override\!bshrinkage{#4}\!dimenG
  \!override\!tshrinkage{#5}\!dimenH
  \ignorespaces}
\def\!hgetshrinkages#1<#2,#3,#4,#5>#6\!nil{%
  \!override\!lshrinkage{#2}\!dimenG
  \!override\!rshrinkage{#3}\!dimenH
  \!override\!bshrinkage{#4}\!dimenE
  \!override\!tshrinkage{#5}\!dimenF
  \ignorespaces}

\def\!shadecolumn{%
  \!dxpos=\!xpos
  \advance\!dxpos -\!xS
  \!removept\!dxpos\!dx
  \!getylimits
  \advance\!ytpos -\!dimenH
  \advance\!ybpos \!dimenG
  \!yloc=\!!yshade
  \ifodd\!parity 
     \advance\!yloc \!dshade
  \fi
  \!lattice\!yloc{2\!dshade}\!ybpos%
    \!countA\!ypos
  \!dimenA=-\!shadexorigin \advance \!dimenA \!xpos
  \loop\!not{\ifdim\!ypos>\!ytpos}
    \!setshadelocation
    \!rotateaboutpivot\!xloc\!yloc%
    \!dimenA=-\!shadexorigin \advance \!dimenA \!xloc
    \!dimenB=-\!shadeyorigin \advance \!dimenB \!yloc
    \kern\!dimenA \raise\!dimenB\copy\!shadesymbol \kern-\!dimenA
    \advance\!ypos 2\!dshade
  \repeat
  \ignorespaces}
 
\def\!qgetylimits{%
  \!dimenA=\!dx\!ytC              
  \advance\!dimenA \!ytB
  \!ytpos=\!dx\!dimenA
  \advance\!ytpos \!ytS
  \!dimenA=\!dx\!ybC              
  \advance\!dimenA \!ybB
  \!ybpos=\!dx\!dimenA
  \advance\!ybpos \!ybS}
 
\def\!lgetylimits{%
  \!ytpos=\!dx\!ytB
  \advance\!ytpos \!ytS
  \!ybpos=\!dx\!ybB
  \advance\!ybpos \!ybS}
 
\def\!vsetshadelocation{
  \!xloc=\!xpos
  \!yloc=\!ypos}
\def\!hsetshadelocation{
  \!xloc=\!ypos
  \!yloc=\!xpos}





\def\!axisticks {%
  \def\!nextkeyword##1 {%
    \expandafter\ifx\csname !ticks##1\endcsname \relax
      \def\!next{\!fixkeyword{##1}}%
    \else
      \def\!next{\csname !ticks##1\endcsname}%
    \fi
    \!next}%
  \!axissetup
    \def\!axissetup{\relax}%
  \edef\!ticksinoutsign{\!ticksinoutSign}%
  \!ticklength=\longticklength
  \!tickwidth=\linethickness
  \!gridlinestatus
  \!setticktransform
  \!maketick
  \!tickcase=0
  \def\!LTlist{}%
  \!nextkeyword}

\def\ticksout{%
  \def\!ticksinoutSign{+}}

\ticksout

\def\nogridlines{%
  \def\!gridlinestatus{\!gridlinestoofalse}}
\nogridlines

\def\loggedticks{%
  \def\!setticktransform{\let\!ticktransform=\!logten}}
\def\unloggedticks{%
  \def\!setticktransform{\let\!ticktransform=\!donothing}}
\def\!donothing#1#2{\def#2{#1}}
\unloggedticks

\expandafter\def\csname !ticks/\endcsname{%
  \!not {\ifx \!LTlist\empty}
    \!placetickvalues
  \fi
  \def\!tickvalueslist{}%
  \def\!LTlist{}%
  \expandafter\csname !axis/\endcsname}

\def\!maketick{%
  \setbox\!boxA=\hbox{%
    \beginpicture
      \!setdimenmode
      \setcoordinatesystem point at {\!zpt} {\!zpt}   
      \linethickness=\!tickwidth
      \ifdim\!ticklength>\!zpt
        \putrule from {\!zpt} {\!zpt} to
          {\!ticksinoutsign\!tickxsign\!ticklength}
          {\!ticksinoutsign\!tickysign\!ticklength}
      \fi
      \if!gridlinestoo
        \putrule from {\!zpt} {\!zpt} to
          {-\!tickxsign\!xaxislength} {-\!tickysign\!yaxislength}
      \fi
    \endpicturesave <\!Xsave,\!Ysave>}%
    \wd\!boxA=\!zpt}
  
\def\!ticksin{%
  \def\!ticksinoutsign{-}%
  \!maketick
  \!nextkeyword}

\def\!ticksout{%
  \def\!ticksinoutsign{+}%
  \!maketick
  \!nextkeyword}

\def\!tickslength<#1> {%
  \!ticklength=#1\relax
  \!maketick
  \!nextkeyword}

\def\!tickslong{%
  \!tickslength<\longticklength> }

\def\!ticksshort{%
  \!tickslength<\shortticklength> }

\def\!tickswidth<#1> {%
  \!tickwidth=#1\relax
  \!maketick
  \!nextkeyword}

\def\!ticksandacross{%
  \!gridlinestootrue
  \!maketick
  \!nextkeyword}

\def\!ticksbutnotacross{%
  \!gridlinestoofalse
  \!maketick
  \!nextkeyword}

\def\!tickslogged{%
  \let\!ticktransform=\!logten
  \!nextkeyword}

\def\!ticksunlogged{%
  \let\!ticktransform=\!donothing
  \!nextkeyword}

\def\!ticksunlabeled{%
  \!tickcase=0
  \!nextkeyword}

\def\!ticksnumbered{%
  \!tickcase=1
  \!nextkeyword}

\def\!tickswithvalues#1/ {%
  \edef\!tickvalueslist{#1! /}%
  \!tickcase=2
  \!nextkeyword}

\def\!ticksquantity#1 {%
  \ifnum #1>1
    \!updatetickoffset
    \!countA=#1\relax
    \advance \!countA -1
    \!ticklocationincr=\!axisLength
      \divide \!ticklocationincr \!countA
    \!ticklocation=\!axisstart
    \loop \!not{\ifdim \!ticklocation>\!axisend}
      \!placetick\!ticklocation
      \ifcase\!tickcase
          \relax 
        \or
          \relax 
        \or
          \expandafter\!gettickvaluefrom\!tickvalueslist
          \edef\!tickfield{{\the\!ticklocation}{\!value}}%
          \expandafter\!listaddon\expandafter{\!tickfield}\!LTlist%
      \fi
      \advance \!ticklocation \!ticklocationincr
    \repeat
  \fi
  \!nextkeyword}

\def\!ticksat#1 {%
  \!updatetickoffset
  \edef\!Loc{#1}%
  \if /\!Loc
    \def\next{\!nextkeyword}%
  \else
    \!ticksincommon
    \def\next{\!ticksat}%
  \fi
  \next}    
      
\def\!ticksfrom#1 to #2 by #3 {%
  \!updatetickoffset
  \edef\!arg{#3}%
  \expandafter\!separate\!arg\!nil
  \!scalefactor=1
  \expandafter\!countfigures\!arg/
  \edef\!arg{#1}%
  \!scaleup\!arg by\!scalefactor to\!countE
  \edef\!arg{#2}%
  \!scaleup\!arg by\!scalefactor to\!countF
  \edef\!arg{#3}%
  \!scaleup\!arg by\!scalefactor to\!countG
  \loop \!not{\ifnum\!countE>\!countF}
    \ifnum\!scalefactor=1
      \edef\!Loc{\the\!countE}%
    \else
      \!scaledown\!countE by\!scalefactor to\!Loc
    \fi
    \!ticksincommon
    \advance \!countE \!countG
  \repeat
  \!nextkeyword}

\def\!updatetickoffset{%
  \!dimenA=\!ticksinoutsign\!ticklength
  \ifdim \!dimenA>\!offset
    \!offset=\!dimenA
  \fi}

\def\!placetick#1{%
  \if!xswitch
    \!xpos=#1\relax
    \!ypos=\!axisylevel
  \else
    \!xpos=\!axisxlevel
    \!ypos=#1\relax
  \fi
  \advance\!xpos \!Xsave
  \advance\!ypos \!Ysave
  \kern\!xpos\raise\!ypos\copy\!boxA\kern-\!xpos
  \ignorespaces}

\def\!gettickvaluefrom#1 #2 /{%
  \edef\!value{#1}%
  \edef\!tickvalueslist{#2 /}%
  \ifx \!tickvalueslist\!endtickvaluelist
    \!tickcase=0
  \fi}
\def\!endtickvaluelist{! /}

\def\!ticksincommon{%
  \!ticktransform\!Loc\!t
  \!ticklocation=\!t\!!unit
  \advance\!ticklocation -\!!origin
  \!placetick\!ticklocation
  \ifcase\!tickcase
    \relax 
  \or 
    \ifdim\!ticklocation<-\!!origin
      \edef\!Loc{$\!Loc$}%
    \fi
    \edef\!tickfield{{\the\!ticklocation}{\!Loc}}%
    \expandafter\!listaddon\expandafter{\!tickfield}\!LTlist%
  \or 
    \expandafter\!gettickvaluefrom\!tickvalueslist
    \edef\!tickfield{{\the\!ticklocation}{\!value}}%
    \expandafter\!listaddon\expandafter{\!tickfield}\!LTlist%
  \fi}

\def\!separate#1\!nil{%
  \!ifnextchar{-}{\!!separate}{\!!!separate}#1\!nil}
\def\!!separate-#1\!nil{%
  \def\!sign{-}%
  \!!!!separate#1..\!nil}
\def\!!!separate#1\!nil{%
  \def\!sign{+}%
  \!!!!separate#1..\!nil}
\def\!!!!separate#1.#2.#3\!nil{%
  \def\!arg{#1}%
  \ifx\!arg\!empty
    \!countA=0
  \else
    \!countA=\!arg
  \fi
  \def\!arg{#2}%
  \ifx\!arg\!empty
    \!countB=0
  \else
    \!countB=\!arg
  \fi}
 
\def\!countfigures#1{%
  \if #1/%
    \def\!next{\ignorespaces}%
  \else
    \multiply\!scalefactor 10
    \def\!next{\!countfigures}%
  \fi
  \!next}

\def\!scaleup#1by#2to#3{%
  \expandafter\!separate#1\!nil
  \multiply\!countA #2\relax
  \advance\!countA \!countB
  \if -\!sign
    \!countA=-\!countA
  \fi
  #3=\!countA
  \ignorespaces}

\def\!scaledown#1by#2to#3{%
  \!countA=#1\relax
  \ifnum \!countA<0 
    \def\!sign{-}
    \!countA=-\!countA
  \else
    \def\!sign{}%
  \fi
  \!countB=\!countA
  \divide\!countB #2\relax
  \!countC=\!countB
    \multiply\!countC #2\relax
  \advance \!countA -\!countC
  \edef#3{\!sign\the\!countB.}
  \!countC=\!countA 
  \ifnum\!countC=0 
    \!countC=1
  \fi
  \multiply\!countC 10
  \!loop \ifnum #2>\!countC
    \edef#3{#3\!zero}%
    \multiply\!countC 10
  \repeat
  \edef#3{#3\the\!countA}
  \ignorespaces}

\def\!placetickvalues{%
  \advance\!offset \tickstovaluesleading
  \if!xswitch
    \setbox\!boxA=\hbox{%
      \def\\##1##2{%
        \!dimenput {##2} [B] (##1,\!axisylevel)}%
      \beginpicture 
        \!LTlist
      \endpicturesave <\!Xsave,\!Ysave>}%
    \!dimenA=\!axisylevel
      \advance\!dimenA -\!Ysave
      \advance\!dimenA \!tickysign\!offset
      \if -\!tickysign
        \advance\!dimenA -\ht\!boxA
      \else
        \advance\!dimenA  \dp\!boxA
      \fi
    \advance\!offset \ht\!boxA 
      \advance\!offset \dp\!boxA
    \!dimenput {\box\!boxA} [Bl] <\!Xsave,\!Ysave> (\!zpt,\!dimenA)
  \else
    \setbox\!boxA=\hbox{%
      \def\\##1##2{%
        \!dimenput {##2} [r] (\!axisxlevel,##1)}%
      \beginpicture 
        \!LTlist
      \endpicturesave <\!Xsave,\!Ysave>}%
    \!dimenA=\!axisxlevel
      \advance\!dimenA -\!Xsave
      \advance\!dimenA \!tickxsign\!offset
      \if -\!tickxsign
        \advance\!dimenA -\wd\!boxA
      \fi
    \advance\!offset \wd\!boxA
    \!dimenput {\box\!boxA} [Bl] <\!Xsave,\!Ysave> (\!dimenA,\!zpt)
  \fi}

\normalgraphs
\catcode`!=12 
%
%
%
\ifx\djatexLoaded\relax\endinput\else\let\djatexLoaded=\relax\fi

%

\newfam\calfont
\font\tencal=eusm10 \font\sevencal=eusm7 \font\fivecal=eusm5
\textfont\calfont=\tencal 
\scriptfont\calfont=\sevencal 
\scriptscriptfont\calfont=\fivecal
\def\cal{\fam=\calfont}	
\newfam\frakturfont
\font\tenfrak=eufm10 \font\sevenfrak=eufm7 \font\fivefrak=eufm5
\textfont\frakturfont=\tenfrak 
\scriptfont\frakturfont=\sevenfrak 
\scriptscriptfont\frakturfont=\fivefrak

\newfam\bbbfont
\font\tenbbb=msbm10 \font\sevenbbb=msbm7 \font\fivebbb=msbm5	
\textfont\bbbfont=\tenbbb 
\scriptfont\bbbfont=\sevenbbb 
\scriptscriptfont\bbbfont=\fivebbb
\def\bbb{\fam=\bbbfont}
\scriptfont\ttfam=\tentt
\scriptscriptfont\ttfam=\tentt

%
\long\def\ifndef#1{\expandafter\ifx\csname#1\endcsname\relax}

\def\mywarning#1{\ifnum\warningsoff=0
	\immediate\write16{l.\the\inputlineno: #1}\fi}
\ifndef{warningsoff}\def\warningsoff{0}\fi

%
\long\def\author#1{\def\Zauthor{#1}}
\long\def\title#1{\def\Ztitle{#1}}

\long\def\date#1{\def\Zdate{#1}}

\def\abstract{\bigbreak\noindent{\bf Abstract}\smallskip\noindent} 

%
\outer\def\section#1\par{\bigbreak\noindent{\bf #1}\nobreak\medskip\noindent}
\outer\def\beginproclaim#1. {\medbreak
	\noindent{\bf #1.\enspace}\begingroup\sl}
\def\endproclaim{\endgroup\par\ifdim\lastskip<\medskipamount 
	\removelastskip\penalty55\medskip\fi}
\def\beginproof#1{\smallskip{\it#1\/}}
\def\endproof{\leavevmode\vrule height0pt width0pt
	depth0pt\nobreak\hfill\proofbox\smallskip}

\def\proofbox{\drawbox{1.2ex}{1.2ex}{.1ex}}
\def\remark#1#2\par{\ifdim\lastskip<\smallskipamount
	\removelastskip\penalty55\smallskip\fi
	{\it #1\/}#2\smallbreak}

%

\def\symbolictags{\def\userawtags{0}}
\ifndef{userawtags}\def\userawtags{1}\fi
\def\tag#1{\ifndef{#1}\mywarning{tag `#1' 
	undefined.}{\hbox{\bf?`}\tt #1\hbox{\bf ?}}\else
	\ifnum\userawtags=0 \csname #1\endcsname
	\else\csname #1\endcsname{ \tt [#1]}\fi\fi}
\def\Tag#1{\tag{#1}}
\def\eqtag#1{(\tag{#1})}
\def\eqTag#1{\ifnum\userawtags=0
	\eqtag{#1}\else
	\lower12pt\hbox{\eqtag{#1}}\fi}
\def\cite#1{[\tag{#1}]}
\def\ecite#1#2{[\tag{#1}, #2]}
\def\nocite#1{}
\def\deftag#1#2{\ifndef{#2}\else\mywarning{tag `#2' defined
	more than once.}\fi 
	\expandafter\def\csname #2\endcsname{#1}}
\def\defcite#1#2{\deftag{#1}{#2}}
\def\bibitem#1{\ifnum\userawtags=0
	\item{[{\tag{#1}}]}\else
	\noindent\hangindent=\parindent\hangafter=1
			\cite{#1}\enskip\fi}

%
\def\rom#1{({\it\romannumeral#1\/})}
\def\Rom#1{{\rm\uppercase\expandafter{\romannumeral#1}}}
\def\a{\alpha}	\def\b{\beta}	\def\c{\gamma}
\def\d{\delta}  \def\e{\varepsilon}	
  \def\cals{{\cal S}}

%
\def\Z{{\bbb Z}} 
\def\R{{\bbb R}} 
\def\C{{\bbb C}} 
\let\sset=\subseteq		%

\let\tensor=\otimes	

\let\isomorphism=\cong		
\let\to=\rightarrow

\def\({\left(} \def\){\right)}
\def\[{\left[} \def\]{\right]}

%

\def\frac#1/#2{\leavevmode\kern.1em\raise.5ex\hbox{\the
	\scriptfont0 #1}\kern-.1em/\kern-.15em\lower.25ex
	\hbox{\the\scriptfont0 #2}} 

\def\drawbox#1#2#3{\leavevmode\vbox{\hrule height #3%
	\hbox{\vrule width #3 height #2\kern #1%
	\vrule width #3}\hrule height #3}}


\def\mathllap#1{\mathchoice
{\llap{$\displaystyle #1$}}%
{\llap{$\textstyle #1$}}%
{\llap{$\scriptstyle #1$}}%
{\llap{$\scriptscriptstyle #1$}}}
\def\set#1#2{\left\{\,#1\mathllap{\phantom{#2}}\mathrel{}\right|\left.#2\mathllap{\phantom{#1}}\,\right\}}

\parskip=0pt

%
%
\magnification=1100
%
%
%
\symbolictags
%

%
%
\deftag{1}{sec-intro}
\deftag{2}{sec-background}
\deftag{2.1}{thm-cartan-hadamard}
\deftag{2.2}{thm-alexandrov}
\deftag{3}{sec-br-covers}
\deftag{3.1}{thm-br-cover-curvature}
\deftag{3.2}{thm-br-cover-basics}
\deftag{3.3}{thm-geodesics-converge}
\deftag{3.4}{thm-if-one-edge-in-branch-locus}
\deftag{3.1}{eq-comparison-of-comparison-complexes}
\deftag{3.2}{eq-def-of-kn}
\deftag{3.3}{eq-bound-on-sum-of-dyadic-steps}
\deftag{3.4}{eq-def-of-lminus}
\deftag{3.5}{thm-geodesics-are-unique}
\deftag{3.1}{fig-cases}
\deftag{3.6}{thm-local-br-cover}
\deftag{4}{sec-iterated}
\deftag{4.1}{thm-boring-lemma}
\deftag{4.2}{thm-paths-with-good-interiors}
\deftag{4.3}{thm-iterated-br-cover}
\deftag{5}{sec-apps}
\deftag{5.1}{thm-hyperplane-complements-aspherical}
\deftag{5.2}{thm-cubic-moduli-space}
\deftag{5.3}{thm-enriques-moduli-space}
\deftag{5.4}{thm-neighborhoods-described}
\deftag{5.5}{thm-weak-homotopy-equiv}
\defcite{1}{alexandrov51:theorem_on_triangles}
\defcite{2}{allcock:ch4-cubic-moduli}
\defcite{3}{allcock:period-lattice-for-enriques-surfaces}
\defcite{4}{bridson95:metric_spaces_of_nonpositive_curvature}
\defcite{5}{charney-davis:branched-covers-of-riemannian-manifolds}
\defcite{6}{davis91:hyperbolization_of_polyhedra}
\defcite{7}{gp-thy-geom-viewpoint}
\defcite{8}{gromov:hyperbolic-groups}
\defcite{9}{horikawa:periods-of-enriques-surfaces-I}
\defcite{10}{januszkiewicz:hyperbolizations}
\defcite{11}{namikawa:periods-of-enriques-surfaces}
\defcite{12}{paulin:hyperbolic_groups_via_hyperbolizing_polyhedra}
\defcite{13}{shah:projective-degeneration-of-enriques-surfaces}
\defcite{14}{gromov-thurston:pinching-constants-for-hyperbolic-manifolds}
\defcite{15}{troyanov:spaces-of-negative-curvature}

%

\def\cp{\C P}		  
\def\ch{\C H}	  	  
\def\moduli{{\cal M}} 	  
\def\U{{\rm U}}	          

\def\x{\kappa}		  
\def\M{M_{\x}^2}	  
\def\sx{\sqrt\x}	  

\def\Xhat{{\widehat X}}	  
\def\Yhat{{\widehat Y}}	  
\def\D{\Delta}		  
\def\Dtil{\tilde{\D}}	  
\def\cat{CAT($\x$)}	  
\def\diam{\mathop{\rm Diam}\nolimits}  
\def\G{\Gamma}		  
\def\btil{\tilde{\b}}	  
\def\dtil{\tilde{\d}}	  
\def\ytil{\tilde y}	  
\def\geo#1{\overline{#1}} 
\def\tri#1{\triangle#1}	  
\def\Util{\tilde{U}}	  
\def\cA{({\bf A})}	  
\def\cB{({\bf B})}	  
\def\cC{({\bf C})}	  
\def\cD{({\bf D})}	  
\def\cE{({\bf E})}	  
\def\cF{({\bf F})}	  
\def\cG{({\bf G})}	  

\def\Tbar{\bar{T}}
\def\Ubar{\bar{U}}
\def\Kbar{\bar{K}}
\def\Bbar{\bar{B}}
\def\abar{\bar{a}}
\def\Abar{\bar{A}}
\def\atil{\tilde{\a}}
\def\xbar{\bar{x}}

\def\S{{\cal S}}  
\def\H{{\cal H}}  
\def\Htil{{\widetilde{\cal H}}}  
\def\s{\sigma}    
\def\Mhat{{\hat M}} 
\def\Nhat{{\hat N}} 
\def\pihat{{\hat\pi}} 
\def\xtil{{\tilde x}} 
\def\ytil{{\tilde y}} 
\def\ztil{{\tilde z}} 
\def\ctil{{\tilde\gamma}} 
\def\Vtil{{\widetilde V}} 
\def\Btil{{\widetilde B}} 
\def\Atil{{\widetilde A}} 
\def\Gammatil{{\widetilde\Gamma}} 

\def\w{\omega} 
\def\E{{\cal E}} 
\def\I{{\cal I}} 
\def\DD{{\cal D}} 

\def\image{\mathop{\rm image}\nolimits}

\noindent
{\bf Metric curvature of infinite branched covers}
\bigskip\noindent
Daniel Allcock\footnote*{Supported in part by an NSF
Postdoctoral Fellowship}
\hfil\break
25 May 1999
\footnote{}{MSC: 53C23 (14J28, 57N65)}
\footnote{}{Keywords: branched cover, ramified cover, Alexandrov
space, cubic surface, Enriques surface}
\bigskip\bigskip

\abstract
We study branched covering spaces in several contexts, proving that
under suitable circumstances the cover satisfies the same upper
curvature bounds as the base space. The first context is 
of a branched cover of an arbitrary metric space that satisfies
Alexandrov's curvature condition \cat, over an arbitrary
complete convex subset. The second context is of a certain sort
of branched cover of a Riemannian manifold over a family of mutually
orthogonal submanifolds. In neither setting do we require that the
branching be locally finite. We apply our results to hyperplane
complements in several complex manifolds of nonpositive
sectional curvature. This implies that two moduli spaces arising
in algebraic geometry are aspherical, namely that of
the smooth cubic surfaces in $\cp^3$ and that of the 
smooth complex Enriques surfaces.

\section{\Tag{sec-intro}. Introduction}

The purpose of this paper is to establish a basic result
in the theory of metric space curvature in the sense of
Alexandrov, together with several applications in algebraic
geometry. A commonly observed phenomenon is that ``taking a
branched cover of almost anything can only introduce negative
curvature''. One can see this phenomenon in elementary examples
using Riemann surfaces, and the idea also plays a role in the
construction \cite{gromov-thurston:pinching-constants-for-hyperbolic-manifolds} of exotic manifolds with negative sectional
curvature. In this paper we work in the maximal generality in
which sectional curvature bounds make sense, namely in the
comparison geometry of Alexandrov. In this setting we will
establish a very strong  theorem concerning the
persistence of upper curvature bounds in branched covers. We include
examples showing that an important completeness hypothesis
cannot be dropped; our examples also disprove several
claims in the literature.

A simple way to build a cover $\Yhat$ of a space $\Xhat$
branched over $\D\sset\Xhat$ is to take any covering space $Y$
of $\Xhat-\D$ and define $\Yhat=Y\cup\D$. We call $\Yhat$ a
simple branched cover of $\Xhat$ over $\D$. Our main result
(theorem~\tag{thm-br-cover-curvature}) states that if $\Xhat$
satisfies Alexandrov's
\cat\ condition and $\D$ is complete and convex then the natural
metric on $\Yhat$ also satisfies \cat. (When $\x>0$ we impose a
minor hypothesis on the diameters of $\Xhat$ and $\Yhat$.) See
section~\tag{sec-background} for a discussion of Alexandrov's
criterion and other background; we follow the conventions of the
book \cite{bridson95:metric_spaces_of_nonpositive_curvature} by
Bridson and Haefliger. Most of section~\tag{sec-br-covers} is
devoted to establishing this theorem. Only partial results can
be obtained without the completeness hypothesis, and we give
these results together with counterexamples when completeness is
not assumed. We also give a local version,
theorem~\tag{thm-local-br-cover}, which allows
one to work with branched covers more complicated than the
simple sort introduced above, and also avoids any diameter
constraints on $\Xhat$ and $\Yhat$. One interesting
twist is that one must take $\D$ to be locally complete in
order to obtain even local results.

The question which motivated this investigation is whether the
moduli space of smooth cubic surfaces in $\cp^3$ is
aspherical (i.e., has contractible universal cover). The answer
is yes, and our argument also establishes the analogous result for the moduli space
of smooth complex Enriques surfaces. To prove these claims, we use the
fact that each of these moduli spaces is known to be covered by a Hermitian
symmetric space with nonpositive sectional curvature, minus an
arrangement of complex hyperplanes. In each case the
hyperplanes have the property that any two of them are
orthogonal wherever they meet. In section~\tag{sec-apps} we show
that such a hyperplane complement is aspherical. We actually
prove a more general result, in the setting of a complete simply
connected Riemannian manifold $\Mhat$ of non-positive sectional
curvature, minus the union $\H$ of suitable submanifolds which
are mutually
orthogonal, complete, and  totally geodesic.

The basic idea is to try to apply standard nonpositive curvature
techniques like the Cartan-Hadamard theorem to the universal
cover $N$ of $M=\Mhat-\H$. The fundamental obstruction is that
$N$ is not metrically complete. This problem can be circumvented
by passing to its metric completion $\Nhat$, but this introduces
problems of its own. First there is the issue of how $N$ and
$\Nhat$ are related. We resolve this by a simple trick that
shows that the inclusion $N\to\Nhat$ is a homotopy
equivalence. The second and more central problem is that $\Nhat$
is not a manifold and not even locally compact. In particular,
one cannot use the techniques of Riemannian geometry. But it is
still a metric space and it turns out to have curvature~$\leq0$,
in the sense that it satisfies Alexandrov's CAT(0) condition
locally. It is then an easy matter to show that $N$ and $\Nhat$
are contractible. In summary, to study the topology of the
manifold $N$ it turns out to be natural and useful to study the
non-manifold $\Nhat$ and use metric-space curvature rather than
Riemannian curvature.

The connection between the very general treatment of metric
curvature and the applications lies in our
study of the curvature of $\Nhat$. For this it suffices to
work locally; the reader should imagine a closed ball $B$ in 
$\C^n$, equipped with
some Riemannian metric, minus the coordinate
hyperplanes. The metric completion of the universal cover of the
hyperplane complement can be obtained by first taking a simple
branched cover of $B$ over one hyperplane, then taking a
simple branched cover of this branched cover over (the preimage
of) the second hyperplane, and so on. If the hyperplanes are
mutually orthogonal and totally geodesic then our main theorem
may be used inductively to study the curvature of the iterated
branched cover. There
are some minor technical issues, which we chase down in section~\tag{sec-iterated}. Note that the base space in each of
the sequence of branched covers fails to be locally compact
(except in the first step). This means that the inductive argument
actually {\it requires} a theorem treating branched covers of
spaces considerably more general than manifolds.

I would like to thank Jim Carlson and Domingo Toledo for their
interest in this work, and for the collaboration
\cite{allcock:ch4-cubic-moduli} that suggested these problems. I
would also like to thank Richard Borcherds, Misha Kapovich and
Bruce Kleiner for useful conversations. Finally, I am grateful to
Brian Bowditch for pointing out an error in an early version.

\section{\Tag{sec-background}. Background}

Let $(X,d)$ be a metric space. A path in $X$ is a continuous map
from a nonempty compact interval to $X$; its initial
(resp. final) endpoint is the image of the least
(resp. greatest) element of this interval. We sometimes describe
a path as being {\it from} its initial endpoint {\it to} its
final endpoint. When we wish to mention its endpoints but not
worry about which is which, we describe the path as {\it
joining} one endpoint {\it and} or {\it with} the other. When we
speak of a point of a path we mean a point in its image. 
If $\c$
is a path in $X$ with domain $[a,b]$ then we define its length
to be
$$ 
\ell(\c)=\sup\set
{\sum_{i=1}^N d(\c(t_{i-1}),\c(t_i))}
{a=t_0\leq t_1\leq\cdots\leq t_N=b,\quad N\geq1}.
$$
This is an element of $[0,\infty]$. We call $X$ a length space
and $d$ a path metric if for all $x,y\in X$ and all $\e>0$ there is
a path of length $<d(x,y)+\e$ joining $x$ and $y$. All of the
spaces in
this paper are length spaces. An important class of length
spaces is that of connected Riemannian manifolds. Given such a manifold
$M$, one defines the `length' of each piecewise differentiable
path in $M$ as a certain integral. Then one defines the distance
between two points of $M$ to be the infimum of the `lengths'
of such paths joining them. By the machinery above, this metric
assigns a length to {\it every} path in $M$.  Happily for the
terminology this agrees with the  `length' when
the latter is defined. See
\ecite{bridson95:metric_spaces_of_nonpositive_curvature}{I.3.15}
for details.

A path $\c$ is called a geodesic parameterized proportionally to
arclength if there exists $k\geq0$ such that
$d(\c(s),\c(t))=k|s-t|$ for all $s$ and $t$ in the domain of
$\c$. We call $\c$ a geodesic if $k=1$.  We sometimes regard two
geodesics as being the same if they differ only by an isometry
of their domains. For example, we use this convention in
assertions about uniqueness of geodesics in $X$. Similarly, we
will sometimes refer to the image of $\c$, rather than $\c$
itself, as a geodesic. 
Sometimes we will even refer to a path as a geodesic when it is
only a geodesic parameterized proportionally to arclength. We say
that $X$ is a geodesic space if any two of its points are joined
by a geodesic. Most of the spaces in this paper are geodesic.  A
subset $Y$ of $X$ is called convex (in $X$) if any two points of
$Y$ are joined by a geodesic of $X$ and every such geodesic
actually lies in $Y$.

A triangle $T$ in $X$ is a triple $(\c_1,\c_2,\c_3)$ of
geodesics of $X$, called the edges of $T$, such that for each
$i$, the final endpoint of $\c_i$ is the initial endpoint of
$\c_{i+1}$; the vertex of $T$ opposite $\c_i$ is defined to be the common
final endpoint of $\c_{i+1}$ and initial endpoint of
$\c_{i-1}$. Here, subscripts should be read
modulo 3. It is possible for a vertex to be opposite more than
one edge; this occurs when an edge of $T$ has length $0$. An
altitude of $T$ is a geodesic of $X$ joining a vertex of $T$ and
a point of an edge opposite it. This terminology does not reduce
to the usual notion of an altitude of a triangle in the
Euclidean plane when $X=\R^2$. Since we will not use the
classical meaning of the term this should cause no confusion.

Now we define the notion of a metric space satisfying a bound on
its curvature. This elegant idea of Alexandrov
\cite{alexandrov51:theorem_on_triangles} captures much of the
flavor of an upper bound on the sectional curvature of a
Riemannian manifold, in the setting of much more general metric
spaces. The idea is that triangles should be thinner than
comparable triangles in some standard space like the 
Euclidean plane.  For each $\x\in\R$,
let $\M$ be the (unique up to isometry) complete simply
connected Riemannian 2-manifold with constant curvature
$\x$. For $\x=0$ or $\x>0$ this space is $\R^2$ or the sphere of
radius $1/\sx$. For $\x<0$ it is the hyperbolic plane equipped with
a suitable multiple of its standard metric. If $T$ is a triangle
in $X$ then a comparison triangle $T'$ for $T$ in $\M$ is a
triangle $(\c_1',\c_2',\c_3')$ in $\M$ such that the domains of
$\c_i$ and $\c_i'$ coincide for each $i$. In particular we
have $\ell(\c_i)=\ell(\c_i')$. Comparison triangles exist unless
$\x>0$ and $T$ has perimeter $>2\pi/\sx$. When they exist they
are unique up to isometry unless
$\x>0$ and $T$ has an edge of length $\pi/\sx$. 
We will arrange things later so that we will not need to worry
about the existence or uniqueness of comparison triangles.
We will
follow the usual convention of taking $2\pi/\sx$ and similar
expressions to represent $\infty$ when $\x\leq0$. This allows
many assertions to be phrased more uniformly.

For each $i$, we say that $\c_i'$ is the edge of $T'$
corresponding to $\c_i$. If  $p$ is a
point of $\c_i$ then the point $p'$ associated to $p$ on the edge
$\c_i'$ is $\c_i'(t)$, where $t$ is such that $\c_i(t)=p$. Note that a
choice of edge containing $p$ is essential for this
construction, since $p$ may lie on more than one edge of $T$. We say
that $T$ satisfies \cat\ if $T$ has perimeter $<2\pi/\sx$ and
for any two edges $\a$ and $\b$ of $T$ and points $p$ on $\a$
and $q$ on $\b$, we have $d(p,q)\leq d(p',q')$. Here $p'$ and
$q'$ are the points of $\a'$ and $\b'$ corresponding to $p$ and
$q$ and $\a'$ and $\b'$ are the edges corresponding to $\a$ and
$\b$ in a comparison triangle for $T$. We say that $X$ satisfies
(or is) \cat\ if $X$ is geodesic and every triangle in $X$ of
perimeter $<2\pi/\sx$ satisfies \cat. The intuitive meaning of
this condition is that $X$ is ``at least as negatively curved''
as $\M$.

We say that $X$ is locally \cat, or has curvature $\leq\x$, if
each point of $X$ has a convex \cat\ neighborhood. Our interest
in spaces with curvature bounded  above stems from the following very
general version of the Cartan-Hadamard theorem, proven by
Bridson and Haefliger
\ecite{bridson95:metric_spaces_of_nonpositive_curvature}{II.5.1}.

\beginproclaim Theorem 
{\Tag{thm-cartan-hadamard}}. 
A complete simply connected length space of curvature
$\leq\x\leq0$ is \cat.
\endproof
\endproclaim 

\noindent
Note that part of the conclusion of
theorem~\tag{thm-cartan-hadamard} is that the space is geodesic,
which is very important and not at all obvious.  This statement
of the theorem implies the version quoted in the introduction:
contractibility follows from the \cat\ condition for $\x\leq0$,
and if the given metric isn't a path metric then it induces one
and the two metrics define the same topology.  
Another very important theorem in Alexandrov's subdivision
lemma, a proof of which appears in \cite{bridson95:metric_spaces_of_nonpositive_curvature}.

\beginproclaim Theorem \Tag{thm-alexandrov} (Alexandrov).
Let $T$ be a triangle in a metric space, with an altitude
$\a$. If $T$ has perimeter $<2\pi/\sqrt\x$ and both of the
triangles into which $A$ subdivides $T$ satisfy \cat, then $T$
also satisfies \cat.
\endproof
\endproclaim

\noindent
Sometimes this is
stated with the conditions that the two subtriangles have
perimeters $<2\pi/\sqrt\x$, but we have made this condition a
part of the definition of \cat.

\section{\Tag{sec-br-covers}. Simple branched covers}

The purpose of this section is to show that under very general
conditions a branched cover satisfies the same upper bounds on
curvature as its base space. Our precise formulation of this
idea is theorem~\tag{thm-br-cover-curvature}. The statement is
slightly stronger than the version given in the introduction
because it turns out that the completeness of the branch locus
is needed for the existence of geodesics but not for the fact
that all triangles in the cover satisfy \cat. The other result
of this section is a local version of this result,
theorem~\tag{thm-local-br-cover}. Since we do not need this
result we will merely state it and give the idea of its proof.

The branched covering spaces we treat here are what
we call simple branched covers. The basic idea is very simple:
one removes a closed subset $\D$ from a length space $\Xhat$,
takes a cover of what is left, and then attaches a copy of $\D$
in the obvious way. Formally, if $\Xhat$ is a length space and
$\D$ is a closed subset of $\Xhat$ then we say that
$\pi:\Yhat\to\Xhat$ is a simple branched cover of $\Xhat$ over
$\D$ if $\Yhat$ is a length space and $\pi$ satisfies the
following two conditions.  First, the restriction of $\pi$ to
$\Yhat-\pi^{-1}(\D)$ must be a locally isometric covering map.
Second, we require that $d(y,z)=d(\pi y,\pi z)$ if at least one
of $y,z\in\Yhat$ lies in $\pi^{-1}(\D)$.  It follows from the
second condition that the restriction of $\pi$ to $\pi^{-1}(\D)$
is an isometry. We will identify $\D$ with its preimage under
$\pi$ and write $X$ and $Y$ for $\Xhat-\D$ and $\Yhat-\D$,
respectively.  It is easy to see that any simple branched cover
is distance non-increasing.

If we
are given $\Xhat$ and $\D$ as above, and
$Y$ is any covering space of $X=\Xhat-\D$, then there is a
unique metric on $\Yhat=Y\cup\D$ such that the obvious map
$\pi:\Yhat\to\Xhat$ is a simple branched cover of $\Xhat$ over
$\D$. This may be constructed as follows. First, each component
of $X$ carries a unique path metric under which its inclusion
into $\Xhat$ is a local isometry. (This uses the fact that $\D$
is closed.) Second, each component of $Y$ carries a natural path
metric, the unique such metric under which the covering map is a
local isometry. Third, for $y,z\in\Yhat$ with at least one of
them in $\D$ we define $d(y,z)=d(\pi y,\pi z)$. Finally, if
$x,z\in Y$ then we define $d(x,z)$ as
$$ 
\inf\Bigl(
\set{d(x,y)+d(y,z)}{y\in\D}
\;\cup\;
\set{\ell(\c)}{\hbox{$\c$ is a path in a component of $Y$ joining
$x$ and $z$}}
\Bigr).
$$
One can check that $d$ is a path metric on $\Yhat$ and that
$\pi$ is a simple branched covering. Our main theorem
is a sufficient condition for $\Yhat$ to be \cat:

\beginproclaim Theorem 
{\Tag{thm-br-cover-curvature}}.  
Suppose $\D$ is a closed convex subset
of a \cat\ space $\Xhat$ and let $\pi:\Yhat\to\Xhat$ be
a simple branched cover of $\Xhat$ over $\D$.  If $\x>0$ then
assume also that $\diam(\Xhat)<\pi/2\sx$ and
$\diam(\Yhat)<2\pi/3\sx$. Then
\item{\rom1}
Every triangle in $\Yhat$ satisfies \cat.
\item{\rom2}
If $\D$ is complete then $\Yhat$ is geodesic and hence \cat.
\endproclaim 

\remark{Example:}
The
completeness condition in \rom2 cannot be dropped, because of the
following example. 
Take $\Xhat$ to be the set of points $(x,y)\in\R^2$ with
$x\geq0$ and $y>0$, together with the point $(1,0)$. Let $\D$ be
the positive $y$-axis. Then $\Xhat$ is a convex subset of
$\R^2$, hence CAT(0), and $\D$ is a closed convex
subset of $\Xhat$. The set $X=\Xhat-\D$ is contractible, so any
cover of it is a union of disjoint copies of it. Taking $Y$ to
be the cover with 2 sheets, $\Yhat$ is isometric to the
upper half plane in $\R^2$ together with the points
$(\pm1,0)$. There is no geodesic joining these two points, so
$\Yhat$ is not a geodesic space. This provides
a counterexample to several assertions in the literature, such
as \ecite{gromov:hyperbolic-groups}{4.3--4.4},
\ecite{januszkiewicz:hyperbolizations}{Lemma~1.1} and 
\ecite{davis91:hyperbolization_of_polyhedra}{Lemma~2.4}.

\remark{Example:}
Although the space $\Yhat$ of the previous example is not
geodesic, it still has curvature $\leq0$. The following example shows
that even this may fail if $\D$ is not complete. We take $\Xhat$
to be the set of points $(x,y,z)$ of $\R^3$ whose first nonzero
coordinate is positive, together with the origin. That is,
$\Xhat$ is the union of an open half-space together with an open
half-plane in its boundary, together with a ray in {\it its}
boundary. We take $\D$ to be the set of points of $\Xhat$ with
vanishing $x$-coordinate, which is the union of the open half-plane
and the ray. Then  $\Xhat$ is a convex
subset of $\R^3$ and $\D$ is closed and convex in
$\Xhat$. As before, any cover of $X=\Xhat-\D$ is a union of copies
of $X$, and we take $Y$ to be the cover with 2 sheets. Then
$\Yhat$ is isometric to the subset of $\R^3$ given by
$$ 
\Yhat=\D\cup\set{(x,y,z)\in\R^3}{x\neq0}, 
$$
equipped with the path metric induced by the Euclidean
metric. It is easy to see that for each $n\geq1$ the points
$(\pm1/n,-1/n,-1/n)$ are joined by no geodesic of
$\Yhat$. Since every neighborhood of $0$ contains such a pair of
points, $0$ has no geodesic neighborhood.

In the proofs below we will use the
following two facts about $\Xhat$. First, geodesics are
characterized by their endpoints. Second, geodesics vary
continuously with respect to their endpoints, by which we mean
that for each $\e>0$ there is a $\d>0$ such that if $d(x,x')<\d$
and $d(y,y')<\d$ for $x,y,x',y'\in\Xhat$ then the geodesic from
$x$ to $y$ is uniformly within $\e$ of the geodesic from $x'$ to
$y'$. These
facts follow from the \cat\ inequalities and the fact that
$\diam(\Xhat)<\pi/\sx$. 
We will ignore all conditions about perimeters of
triangles in $\Yhat$ being less than $2\pi/\sx$ because we have
bounded $\diam(\Yhat)$ in order to guarantee that all
triangles in $\Yhat$ satisfy this condition. We chose the bounds
on $\diam(\Xhat)$ and $\diam(\Yhat)$ out of convenience; one could
probably weaken them, although most theorems about \cat\ spaces
require some sort of extra condition when $\x>0$. Of course, if
$\diam(\Xhat)<\pi/3\sx$ then the condition on $\diam(\Yhat)$
follows automatically.

We begin with some elementary properties of geodesics in
$\Yhat$, and then show that under special circumstances they
vary continuously with respect to their endpoints.

\beginproclaim Lemma 
{\Tag{thm-br-cover-basics}}.  Under the hypotheses of
theorem~\tag{thm-br-cover-curvature}, we have the following:
\item{\rom1}
If $\c$ is a geodesic of $\Xhat$ meeting $\D$ and $w,z\in\Yhat$
lie over the endpoints of $\c$, then there is a unique path in
$\Yhat$ joining $w$ with $z$ and projecting to $\c$.
\item{\rom2}
A path in $\Yhat$ projecting to a geodesic of $\Xhat$ is the
unique geodesic with its endpoints.
\item{\rom3}
If $x\in\Yhat$ and $y\in\D$ then there is a unique geodesic of
$\Yhat$ joining them, and it projects to a geodesic of $\Xhat$.
\item{\rom4}
A geodesic of $\Yhat$ that misses $\D$ projects to a geodesic of $\Xhat$.
\item{\rom5}
If $x\in Y$ then the set of points of $\Yhat$ that may be joined
with $x$ by a geodesic of $\Yhat$ that misses $\D$ is open.
\par
\endproclaim

\beginproof{Proof:}
\rom1
Because $\D$ is convex in $\Xhat$, $\c$ meets $\D$ in an
interval, with endpoints say $x$ and $y$. There is clearly a
unique lift of this interval. If $\pi w\neq x$ then by covering
space theory the half-open segment of $\c$ from $\pi w$ to $x$
has a unique lift to $Y$ beginning at $w$. Similarly, there is a
unique lift of the
segment from $y$ to $\pi z$. It is obvious
that these lifts fit together to form a lift of $\c$.

\rom2
This follows from the uniqueness of geodesics in $\Xhat$ and the
uniqueness of their lifts with specified endpoints in $\Yhat$,
which in turn follows from \rom1 for geodesics that meet $\D$ and
from covering space theory for those that do not.

\rom3
This follows from \rom1 and \rom2 by lifting a geodesic of
$\Xhat$ joining $\pi x$ and $\pi y$.

\rom4
Suppose $\c$ is a geodesic of $\Yhat$ from $x$ to $z$ that
misses $\D$, and that $\pi\c$ is not a geodesic of $\Xhat$. 
Consider the homotopy $\G$ from
$\pi\c$ to the constant path at $\pi x$ given by retraction
along geodesics. Suppose first that $\G$ meets $\D$. Then there
is a point $y$ of $\pi\c$ that is joined with $\pi x$ by a geodesic $\b$
of $\Xhat$ that meets $\D$. By \rom1, there is a lift $\btil$ of $\b$ from $x$
to the point $\ytil$ of $\c$ lying over $y$. By \rom2, $\btil$
is the unique geodesic of $\Yhat$ with these endpoints. But then
the subsegment of $\c$ from $x$ to $\ytil$ must coincide with
$\btil$, contradicting the fact that $\c$ misses $\D$. Now
suppose $\G$ misses $\D$. Consider
the geodesic $\d$ of $\Xhat$ from $\pi x$ to $\pi z$, which
is shorter than $\c$. 
Since $\d$ is the track of $\pi z$, we
may regard $\G$ as a homotopy rel endpoints between $\pi\c$ and
$\d$. This lifts to a homotopy between $\c$ and a lift $\dtil$
of $\d$ that joins $x$ and $z$. Then $\ell(\dtil)<\ell(\c)$,
contradicting the hypothesis that $\c$ is a geodesic.

\rom5
Suppose $\c$ is a geodesic of $\Yhat$ from $x$ to some point $y$
of $\Yhat$, and that $\c$ misses $\D$. By \rom4, $\pi\c$ is a
geodesic of $\Xhat$. Since geodesics of $\Xhat$ depend
continuously upon their endpoints, there is an open ball $U$ of
radius $r>0$ about $\pi y$ such that the geodesics of $\Xhat$
from $\pi x $ to the various points of $U$ are all uniformly
within $d(\image(\c),\D)$ of $\pi\c$. 
By replacing $r$ by a smaller number if necessary, we
may also suppose that the open $r$-ball $\Util$ about $y$ maps
isometrically onto its image. Now, if $y'$ lies in $\Util$
then consider the homotopy along geodesics from $\pi\c$ to the
geodesic $\b$ from $x$ to $\pi y'$.  This misses $\D$, so it
lifts to a homotopy from $\c$ to a lift $\btil$ of $\b$.
Considering the length of the track of $y$, we see that the
final endpoint of $\btil$ lies within $r$ of $y$, so that it
must coincide with $y'$. Finally, $\btil$ is a geodesic by
\rom2. This shows that every point of $\Util$ is joined with $x$
by a geodesic that misses $\D$.
\endproof

\beginproclaim Lemma 
{\Tag{thm-geodesics-converge}}. 
Under the hypotheses of theorem~\tag{thm-br-cover-curvature},
suppose $x_n$ and $y_n$ are sequences in $Y$ converging to
points $x$ and $y$ of $\Yhat$, respectively. Suppose also that
for each $n$ there is a geodesic of $\Yhat$ from $x_n$ to $y_n$
that misses $\D$, and let $\c_n:[0,1]\to\Yhat$ be a parameterization
of this geodesic proportional to arclength. Then there is a
geodesic of\/ $\Yhat$ from $x$ to $y$, and the $\c_n$ converge
uniformly to (the obvious reparameterization of) it.
\endproclaim

\beginproof{Proof:}
By lemma~\tag{thm-br-cover-basics}\rom4, each $\pi\c_n$ is a
geodesic of $\Xhat$. By the continuous dependence of geodesics
in $\Xhat$ on their endpoints, the $\pi\c_n$ converge uniformly
to a (suitably parameterized) geodesic $\b$ from $\pi x$ to $\pi
y$. We distinguish two cases. First, suppose that $\b$ misses
$\D$.  Then there is a unique lift $\btil$ of $\b$ with
$\btil(0)=x$, and this is a (reparameterized) geodesic by
\tag{thm-br-cover-basics}\rom2. We claim that the $\c_n$ converge uniformly to
$\btil$ and that $\btil(1)=y$, so that $\btil$ is the required
geodesic. The second claim follows from the first. To see
convergence, first choose a constant $\d>0$ such that $\d<
d(\image(\b),\D)$ and such that the open ball of radius $\d$
about $x$ maps isometrically onto its image.  By discarding
finitely many terms of the sequence, we may suppose that all the
$x_n$ are within $\d$ of $x$ and that all the $\pi\c_n$ are
uniformly within $\d$ of $\b$. 

We claim that for each $n$, the uniform distance between $\c_n$
and $\btil$ equals the uniform distance between $\pi\c_n$ and
$\b$. This clearly implies that the $\c_n$ converge uniformly to
$\btil$. To see the claim, simply construct the homotopy along
geodesics in $\Xhat$ from $\pi\c_n$ to $\b$ and lift this to a
homotopy from $\c_n$ to some lift of $\b$, and then argue as in
lemma~\tag{thm-br-cover-basics}\rom5 that this lift coincides with $\btil$.

On the other hand, suppose that $\b$ meets $\D$. Then parts
\rom1 and \rom2 of lemma~\tag{thm-br-cover-basics} show that
there is a unique lift $\btil$ of $\b$ with $\btil(0)=x$ and
$\btil(1)=y$, and that this is a (reparameterized) geodesic of
$\Yhat$. We must show that the $\c_n$ converge uniformly to
$\btil$. Let $u$ (resp. $w$) be the least (resp. greatest)
element of $[0,1]$ whose image under $\b$ lies in $\D$.  It is
easy to see that the $\c_n$ converge uniformly to $\btil$ on
$[u,w]$. (One just uses the fact that the distance between an
element of $\Yhat$ and an element of $\D\sset\Yhat$ coincides
with the distance between their projections.) From this and the
fact that the $\c_n$ are all parameterized proportionally to
arclength, one obtains the following slightly stronger
statement: for each $\e>0$ there exist an $N$ and a $\d>0$ such
that for all $n>N$ and all $v\in[u-\d,w+\d]$,
$d(\c_n(v),\btil(v))<\e$. Therefore it suffices to prove that
for all $\d>0$, the $\c_n$ converge uniformly on $[0,u-\d]$ and
on $[w+\d,1]$. Since the images of these intervals under $\b$
are disjoint from $\D$ we can use the same argument as in the
case that $\b$ missed $\D$.
\endproof

\beginproclaim Lemma 
{\Tag{thm-if-one-edge-in-branch-locus}}. 
A triangle of $\Yhat$ with an edge in $\D$ satisfies \cat.
\endproclaim

\beginproof{Proof:}
We write $T$ for the triangle, and represent the edge in $\D$ as
$A:[0,1]\to\D$, a geodesic parameterized proportionally to
arclength. Let $a$ be the vertex of $T$ opposite $A$, and for
$t\in[0,1]$ let $B^t$ be the geodesic from $a$ to $A(t)$,
parameterized proportionally to arclength. In particular, the
edges of $T$ are $A$, $B^0$ and $B^1$. Let $T'$ be a comparison
triangle in $\M$ for $T$, and suppose that $p$ and $q$ are
points of given edges of $T$. Let $p'$ and $q'$ be the
corresponding points of the corresponding edges of $T'$, 
and let $k=d_{T'}(p',q')$. We must
show $d_{\Yhat}(p,q)\leq k$. If one of $p$ and $q$
lies on $A$ then we are done, because
$$ 
d_{\Yhat}(p,q)=d_{\Xhat}(\pi p, \pi q)\leq  k\;,
$$
where we have used the definition of a simple branched cover,
the fact that $\Xhat$ is \cat, and the fact that $T'$ is a
comparison triangle for $\pi T$ as well as for $T$.

Now we consider the case in which neither $p$ nor $q$ is given
as lying on $A$. To avoid trivialities we suppose that they lie
on different edges of $T$, so we may suppose that $p$ lies on
$B^0$ and $q$ on $B^1$. The obvious idea is to construct the
geodesic in $\Xhat$ joining $x_0=\pi p$ with $x_1=\pi q$, and
then lift it to a geodesic of $\Yhat$. The problem is that while
we may always lift the geodesic, there is no guarantee that the
lift will join $p$ and $q$. We will circumvent this problem by
joining $x_0$ to $x_1$ by a path $\a$ that may fail to be a
geodesic, but will have length~$\leq k$. Our path will have
the virtue of lying in the `surface' $S$ swept out by the geodesics
$\pi B^t$, which will allow us to lift it to a path from $p$ to
$q$. By the continuous dependence of geodesics on their
endpoints, $(t,z)\mapsto \pi B^t(z)$ is continuous on
$[0,1]\times[0,1]$, and so $S$ is compact. 

We will need some ``comparison complexes'' $\Kbar_n$ as well as
the comparison triangle $T'$. For $0\leq t\leq u\leq1$ we define
$T(t,u)$ to be the triangle with edges $B^t$, $B^u$ and
$A|_{[t,u]}$. For each $n=0,1,2,\ldots$ we let $D_n$ be the set
of dyadic rational numbers in $[0,1]$ of the form
$t_{n,i}=i/2^n$. For each $i=1,\ldots,2^n$
we define $\Tbar_{n,i}$ to be the comparison triangle in $\M$
for $T(t_{n,i-1},t_{n,i})$. We write $\Bbar_{n,i}^-$ and
$\Bbar_{n,i}^+$ for the edges of $\Tbar_{n,i}$ corresponding to
$B^{t_{n,i-1}}$ and $B^{t_{n,i}}$, $\Abar_{n,i}$ for the edge
corresponding to $A|_{[t_{n,i-1},t_{n,i}]}$, and $\abar_{n,i}$
for the vertex corresponding to $a$. We take $\Ubar_{n,i}$ to be
the convex hull of $\Tbar_{n,i}$ in $\M$. Finally, we define
$\Kbar_n$ as the union of disjoint copies of the $\Ubar_{n,i}$,
subject to the identification of the the segment $\Bbar_{n,i}^+$
in $\Ubar_{n,i}$ with $\Bbar_{n,i+1}^-$ in $\Ubar_{n,i+1}$, in
such a way that $\abar_{n,i}$ is identified with
$\abar_{n,i+1}$, for each $i=1,\ldots,2^n-1$. In short,
$\Kbar_n$ is a `fan' of $2^n$ triangular pieces cut from $\M$,
although some of these pieces may degenerate to segments.  We
equip $\Kbar_n$ with its natural path metric. The paths
$\Abar_{n,i}$ in the $\Ubar_{n,i}$ fit together to form a path
$\Abar_n:[0,1]\to\Kbar_n$.  The vertices $\abar_{n,i}$ are
identified with each other, resulting in a single point
$\abar_n$ of $\Kbar_n$. If $t\in D_n$ then we let
$\Bbar_n^t:[0,1]\to\Kbar_n$ be the geodesic from $\abar_n$ to
$\Abar_n(t)$, parameterized proportionally to arclength. These
paths, together with $\Abar_n$, form the `1-skeleton' of
$\Kbar_n$. 

We claim that $\Abar_n$ is a geodesic of $\Kbar_n$ for all $n$.
Otherwise, a simple application of the \cat\ property of
$\Xhat$ would show that $\pi A$ failed to be a geodesic. We use
this to deduce that points of the `1-skeleton' of $\Kbar_n$ are
at least as far apart as the corresponding points of
$\Kbar_{n+1}$. To make this precise, observe that if $t\in D_n$
then $\pi B^t$ and $\Bbar^t_n$ have the same length. Therefore
to each point $x$ of $\pi B^t$ we may associate a point $\xbar$
of $\Bbar^t_n$ and vice-versa. The relationship is
$d_{\Xhat}(\pi a,x)=d_{\Kbar_n}(\abar_n,\xbar)$. If $t$ also
lies in $D_m$ then we can identify $\Bbar^t_n$ with $\Bbar^t_m$
in a similar way. We claim that if $u$ and $w$ lie in $D_{n-1}$
and $b$ and $c$ are points on $\Bbar^u_{n-1}$ and
$\Bbar^w_{n-1}$, with corresponding points $\b$ and $\c$ on
$\Bbar^u_n$ and $\Bbar^w_n$, then 
$$
d_{\Kbar_n}(\b,\c)\leq d_{\Kbar_{n-1}}(b,c)\;.
\eqno\eqTag{eq-comparison-of-comparison-complexes}
$$
To prove this it suffices to treat the case in which $u$ and $w$
are consecutive elements of $D_{n-1}$. Since $\Abar_n|_{[u,w]}$
is a geodesic, we may consider the geodesic triangle in
$\Kbar_n$ with this edge together with $\Bbar_n^u$ and
$\Bbar_n^w$. This satisfies \cat\ because we may subdivide it
along the altitude $\Bbar_n^v$ (where $v=(u+w)/2$), into two
triangles which are pieces of $\M$ and therefore obviously
\cat. Furthermore, as a comparison triangle we may take the
triangle in $\Kbar_{n-1}$ bounded by $\Abar_{n-1}|_{[u,w]}$,
$\Bbar_{n-1}^u$ and $\Bbar_{n-1}^w$, since this triangle is also
a piece of $\M$. Then
\eqtag{eq-comparison-of-comparison-complexes} follows
immediately.

We write $D=\cup D_n$ for the set of all dyadic rational numbers
in $[0,1]$. We will use the $\Kbar_n$ to construct a point $x_u$
on $\pi B^u$ for each $u\in D$. Then we will string the $x_u$
together to build the path $\a$. We have already defined $x_0=\pi
p$ and $x_1=\pi q$. We will sometimes write $x_{n,i}$ for
$x_{i/2^n}$, so we have just defined $x_{0,0}=x_0$ and
$x_{0,1}=x_1$. Supposing that all the $x_{n-1,i}$ have been
defined, we define the $x_{n,j}$ as follows. We have already
defined the $x_{n,j}$ for even $j$, namely
$x_{n,j}=x_{n-1,j/2}$. If $j$ is odd then take $u=(j-1)/2^n$,
$v=j/2^n$ and $w=(j+1)/2^n$, and consider the points $\xbar_u$
and $\xbar_w$ of $\Kbar_n$ that lie on $\Bbar_n^u$ and
$\Bbar_n^w$ and correspond to $x_u$ and $x_w$. We construct the
geodesic of $\Kbar_n$ joining $\xbar_u$ and $\xbar_w$, and let
$\xbar_v$ be any point of $\Bbar_n^v$ that it meets. Such an
intersection point exists by the construction of
$\Kbar_n$. (Typically, there will be a unique intersection
point, but this can fail if one of $\Tbar_{n,j}$ and
$\Tbar_{n,j+1}$ degenerates to a segment.) We take $x_v$ to be
the point of $\pi B^v$ corresponding to $\xbar_v$.

For each $n$ we write $\xbar_{n,i}$ for the point of $\Bbar^u_n$
corresponding to $x_u$, where $u=i/2^n$. (In particular,
$\xbar_{n,i}$ and $\xbar_{n+1,2i}$ are points of different
spaces, but both correspond to $x_u$.) Consider the sum
$$
k_n=\sum_{i=1}^{2^n}d_{\Kbar_n}(\xbar_{n,i-1}\,,\,\xbar_{n,i})\;.
\eqno\eqTag{eq-def-of-kn}
$$
We claim that $k_n\leq k_{n-1}$. To see this, observe that
$$ 
k_n=\sum_{i=2,4,\ldots,2^n}
d_{\Kbar_n}(\xbar_{n,i-2}\,,\,\xbar_{n,i}) 
$$
because of the construction of the $x_{n,i}$ for odd $i$. By
\eqtag{eq-comparison-of-comparison-complexes}, when $i$ is  even we have
$$ 
d_{\Kbar_n}(\xbar_{n,i-2}\,,\,\xbar_{n,i})
\leq
d_{\Kbar_{n-1}}(\xbar_{n-1,(i-2)/2}\,,\,\xbar_{n-1,i/2})\;, 
$$
and $k_n\leq k_{n-1}$ follows. We immediately obtain $k_n\leq
k_0=k$ for all $n$. By applying the \cat\ inequality for $\Xhat$
we see that for all $n$,
$$ 
\sum_{i=1}^{2^n}d(x_{n,i-1}\,,\,x_{n,i})\leq k\;. 
$$
It follows immediately that if $u_0,\ldots,u_j$ is any
increasing sequence of dyadic rationals, then
$$
d(x_{u_0},x_{u_1})+\cdots+
d(x_{u_{j-1}},x_{u_j})\leq k\;.
\eqno\eqTag{eq-bound-on-sum-of-dyadic-steps}
$$

We are now ready to string the $x_u$'s together into a
path. There is a technical complication, which we will work
around in the next few paragraphs. Specifically, the map
$u\mapsto x_u$ might not be continuous on $D$. This can happen
if some of the comparison triangles $\Tbar_{n,i}$ are
degenerate. However, if $t\in(0,1]$ then as the $u\in
D$ approach $t$ from the left, the $x_u$ do converge to a limit
$L(t)$. Similarly, if $t\in[0,1)$ then as the $u\in D$ approach
$t$ from the right, the $x_u$ converge to a limit $R(t)$. We
will treat $L(t)$; the discussion of $R(t)$ is similar. Certainly
there is some sequence $u_i$ of dyadic rationals approaching $t$
from below, such that the $x_{u_i}$ converge, since the $x_u$
all lie in the compact set $S$. We will call this limit
$L(t)$. Now we show that if $u_j'$ is {\it any} sequence in $D$
approaching $t$ from below then the $x_{u_j'}$ converge to
$L(t)$. For otherwise we could suppose (by passing to a
subsequence) that the $x_{u_j'}$ converge to some other other
limit. Then by interleaving terms of the sequences $u_i$ and
$u_j'$ we could violate
\eqtag{eq-bound-on-sum-of-dyadic-steps}. This establishes the
existence of the left and right limits $L(t)$ and $R(t)$. It is
obvious that $L(t)$ and $R(t)$ lie on $\pi B^t$. For
completeness we define $L(0)=x_0$ and $R(1)=x_1$.

Next, we claim that if $t_1,t_2,\ldots$ is an increasing
sequence in $[0,1]$ with limit $t$, then
$L(t)=\lim_{n\to\infty}R(t_n)$. If this failed then there would
be such a sequence that converged to some point other than
$L(t)$. But then for each $n$ we could choose $u_n\in D$ such
that $t_n<u_n<t_{n+1}$ and $d(x_{u_n},R(t_n))<1/n$. Then the
$x_{u_n}$ would converge to a point other than $L(t)$, while the
$u_n$ approach $t$ from below, a contradiction. A symmetric
argument shows that if $t_1,t_2,\ldots$ is a decreasing sequence
with limit $t$ then $R(t)=\lim_{n\to\infty}L(t_n)$.

Our path $\a$ will pass through all the points $L(t)$ and $R(t)$
in order.  To accomplish this, we define for $t\in(0,1]$ the
quantity
$$
l^-(t)=\sup\left\{\sum_{i=1}^n
d(L(t_{i-1}),R(t_{i-1}))+d(R(t_{i-1}),L(t_i))\right\}\;,
\eqno\eqTag{eq-def-of-lminus}
$$
where the supremum is over all increasing sequences
$0=t_0<\cdots<t_n=t$.  For completeness we define $l^-(0)=0$.
Then for $t\in[0,1]$ we define $l^+(t)=l^-(t)+d(L(t),R(t))$.  To
motivate these definitions, we mention that the length of the
subpath of $\a$ from $x_0$ to $L(t)$ (resp. $R(t)$) will be
$l^-(t)$ (resp. $l^+(t)$).  It is obvious that if $t<t'$ then
$l^-(t)\leq l^+(t)\leq l^-(t')\leq l^+(t')$.  Furthermore,
$l=l^+(1)$ satisfies $l\leq k$. To see this, consider
an increasing sequence $0=t_0<\cdots<t_n=1$ 
such that
$$ 
d(L(t_0),R(t_0))+d(R(t_0),L(t_1))+\cdots+d(L(t_n),R(t_n))
$$
approximates $l$. We may approximate each
$t_i$ (except $t_0$)  by a dyadic rational $u_i$ smaller than $t_i$, and
each $t_i$ (except $t_n$) by a dyadic rational
$v_i$ larger than $t_i$. We may do this in such a way that 
the sequence
$0,v_0,u_1,v_1,\ldots,u_{n-1},v_{n-1},u_n,1$ 
is increasing. Then the $u_i$ approximate the $L(t_i)$ and the
$v_i$ approximate the $R(t_i)$. It follows that
$$ 
d(x_0,x_{v_0})+d(x_{v_0},x_{u_1})+d(x_{u_1},x_{v_1})+\cdots+d(x_{u_n},x_1)
$$
approximates $l$, and then $l\leq k$ follows from
\eqtag{eq-bound-on-sum-of-dyadic-steps}.
Finally, we claim that if $t\in(0,1]$ then
$l^-(t)=\sup_{t'<t}l^+(t')$ and if $t\in[0,1)$ then
$l^+(t)=\inf_{t'>t}l^-(t')$. This follows from the 
relationship between the $L(t)$ and the $R(t)$, together with
the fact that $l$ is finite.

Now we build $\a$. One can check that there is a unique function
$\a:[0,l]\to\Xhat$ such that for each $t\in[0,1]$ the
restriction of $\a$ to $[l^-(t),l^+(t)]$ is the geodesic from
$L(t)$ to $R(t)$. It follows from the relations between the
$L(t)$ and the $R(t)$ that $\a$ is continuous, and from the
definitions of $l^\pm(t)$ that $\a$ is parameterized by
arclength. In particular, $\ell(\a)=l\leq k$. Finally, $\a$ lies
in $S$ since  $L(t)$ and $R(t)$ lie on $\pi B^t$ for each $t$.

Now we will lift $\a$ to $\Yhat$. Suppose first that $\a$ misses
$\D$. If a point $x$  of $\a$  lies on $\pi B^t$, then the
subsegment of $\pi B^t$ from $\pi a$ to $x$ misses $\D$, for
otherwise to convexity of $\D$ would force $x\in \D$. We may
regard the retraction of $\a$ along geodesics to $\pi a$ as a
homotopy rel endpoints between $\a$ and the path $\b$ which
travels along $\pi B^0$ from $x_0$ to $\pi a$ and then along
$\pi B^1$ from $\pi a$ to $x_1$. Of course $\b$ lifts to a path
$\btil$ from $p$ to $q$, and since the homotopy misses $\D$ it
may also be lifted. Therefore there is a lift $\atil$ of $\a$
from $p$ to $q$, with length $l\leq k$, as desired. On the other
hand, if $\a$ meets $\D$ then the lifting is even easier. One
defines $\atil$ on $\atil^{-1}(\D)$ in the obvious way, and then
one defines the rest of $\atil$ by lifting each component of
$\a^{-1}(\Xhat-\D)$ however one desires, subject to the
conditions $\atil(0)=p$ and $\atil(l)=q$.
\endproof 

\beginproclaim Lemma 
\Tag{thm-geodesics-are-unique}. 
At most one geodesic joins any two given points of $\Yhat$.
\endproclaim

\beginproof{Proof:}
Suppose $x,y\in\Yhat$; we claim that there is at most one
geodesic joining them. If either $x$ or $y$ lies in $\D$ then
lemma~\tag{thm-br-cover-basics}\rom3 applies. If they are joined
by a geodesic missing $\D$ then its uniqueness follows from
lemma~\tag{thm-br-cover-basics}\rom4 and \rom2. So it suffices
to consider the case with $x,y\in Y$ such that every geodesic
joining them meets $\D$. Let $\c$ and $\d$ be two such
geodesics, meeting $\D$ in points $c$ and $d$ respectively; we
will prove $\c=\d$. The (unique) geodesic triangle with vertices
$x$, $c$ and $y$ satisfies \cat\ by
lemma~\tag{thm-if-one-edge-in-branch-locus}.  So does the
geodesic triangle with vertices $y$, $c$ and $d$.

We now apply Alexandrov's subdivision lemma
to the `bigon' formed by $\c$ and $\d$. Taking $T$ to
be the triangle with edges $\d$ and the subsegments of $\c$
joining $c$ to each of $x$ and $y$, and the altitude to be the
geodesic joining $c$ and $d$, we see that $T$ satisfies
\cat. Since $\ell(\c)=\ell(\d)$, the comparison triangle
degenerates to a segment, and the \cat\ inequality immediately
implies $\c=\d$.
\endproof

The proofs of the two parts of theorem~\tag{thm-br-cover-curvature} are independent of
each other.

\beginproof{Proof of theorem~\tag{thm-br-cover-curvature}\rom1:}
Cases \cA--\cG\ below show that various sorts of triangles in
$\Yhat$ satisfy \cat. These constitute a proof because every
triangle is treated either by case \cA\ or by case \cG.  In each
case $T$ is a triangle with vertices $A$, $B$ and $C$. For two
points $P$ and $Q$ of $\Yhat$ that are joined by a geodesic we
write $\geo{PQ}$ for the geodesic joining them. For three points
$P$, $Q$ and $R$ of $\Yhat$, any two of which are joined by a
geodesic, we write $\tri{PQR}$ for the geodesic triangle with
edges $\geo{PQ}$, $\geo{QR}$ and $\geo{RP}$.
Figure~\tag{fig-cases} illustrates the arguments for a few of
the cases.  Most of the cases use Alexandrov's lemma. 
Because $\Yhat$ might not be a geodesic space we have to prove
the existence of all the geodesics we introduce, which
complicates the argument. At a first reading one should simply
assume that all needed geodesics exist.

\topinsert
\hfil
\beginpicture
\put {\epsfbox[97 497  486 594]{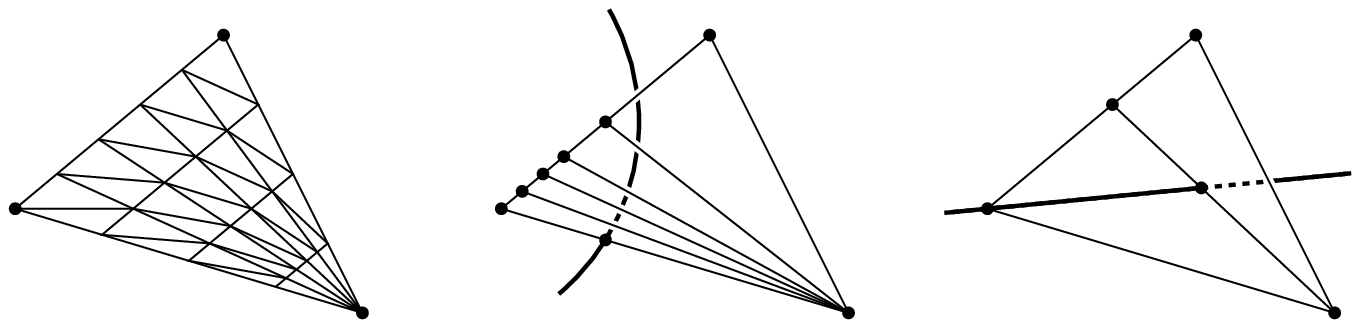}} [bl] at 0 0
\setcoordinatesystem units <1bp,1bp> point at 97 497
\put {$A$} [l] <0pt,0pt> at 201.8 500.0
\put {$B$} [r] <-1pt,0pt> at 98.2 530.0
\put {$C$} [b] <0pt,2pt> at 160.0 581.8
\put {\bf(A)} [c] <0pt,0pt> at 145.0 485.0
\put {$A$} [l] <0pt,0pt> at 341.8 500.0
\put {$B$} [tr] <-1pt,0pt> at 238.2 530.0
\put {$C$} [b] <0pt,2pt> at 300.0 581.8
\put {$B_1$} [br] <0pt,0pt> at 268.2 555.0
\put {$B_n$} [br] <0pt,0pt> at 244.2 535.0
\put {$\D$} [tl] <0pt,0pt> at 256.563 505.467
\put {\bf(B)} [c] <0pt,0pt> at 285.0 485.0
\put {$A$} [l] <0pt,0pt> at 481.8 500.0
\put {$B$} [br] <0pt,2pt> at 378.2 530.0
\put {$C$} [b] <0pt,2pt> at 440.0 581.8
\put {$B'$} [br] <0pt,0pt> at 414.2 560.0
\put {$\D$} [b] <0pt,2pt> at 484.72 540.2
\put {\bf(B)} [c] <-1pt,-3pt> at 440.8 558.0
\put {\bf(C)} [c] <-5pt,3pt> at 419.8 540.5
\put {\bf(C)} [c] <0pt,-2pt> at 435.8 525.5
\put {\bf(E)} [c] <0pt,0pt> at 425.0 485.0
\endpicture
\hfil

\medskip\narrower\narrower\noindent
Figure~\Tag{fig-cases}. Three of the cases in the proof of
theorem~\tag{thm-br-cover-curvature}\rom1. Each picture
represents a triangle in $\Yhat$.  The bold letters in case \cE\
indicate the earlier cases to which the problem is reduced.
\endinsert

\cA\
Suppose no altitude from $A$ meets $\D$. The set of points of
$\geo{BC}$ to which there is a geodesic from $A$ is nonempty
because it contains $C$. Because no altitude from $A$ meets
$\D$, this set is open by
lemma~\tag{thm-br-cover-basics}\rom5. It is also closed
(lemma~\tag{thm-geodesics-converge}), hence all of $\geo{BC}$.
For each $P$ in $\geo{BC}$, let $\c_P:[0,1]\to\Yhat$ be a
parameterization proportional to arclength of the geodesic
$\geo{AP}$ from $A$ to $P$. For each fixed $s\in [0,1]$ the map
$\G:\geo{BC}\times[0,1]\to\Yhat$ given by $(P,s)\mapsto\c_P(s)$
is continuous in $P$, by lemma~\tag{thm-geodesics-converge}. For
each fixed $P$ the map is lipschitz as a function of $s$, with
lipschitz constant $\diam(T)$.  It follows that $\G$ is jointly
continuous in $P$ and $s$. The fact that $T$ satisfies \cat\ now
follows from a standard subdivision argument, like that of
\ecite{paulin:hyperbolic_groups_via_hyperbolizing_polyhedra}{p.~328}.

\cB\
Suppose that the only altitude from $A$ meeting $\D$ is
$\geo{AB}$. If $B=C$ then $T$ degenerates to a segment (by the
uniqueness of geodesics) and therefore automatically satisfies
\cat. If $B\neq C$ then arguing as in the previous case we see
that each point $P$ of $\geo{BC}$ is joined with $A$ by a
geodesic. We choose a sequence of points $B_n$ of
$\geo{BC}-\{B\}$ approaching $B$. For each $n$, $\tri{AB_nC}$
satisfies \cat\ by case \cA. By
lemma~\tag{thm-geodesics-converge}, the geodesics $\geo{AB_n}$
converge uniformly to $\geo{AB}$. As a uniform limit of
triangles that satisfy \cat, $\tri{ABC}$ does also.

\cC\
Suppose $\D$ contains two vertices of $T$. This is
lemma~\tag{thm-if-one-edge-in-branch-locus}.

\cD\
Suppose $\D$ contains a vertex of $T$ and also a point of an
opposite side. Then there is a geodesic joining these points, by
lemma~\tag{thm-br-cover-basics}\rom3. Subdivide $T$ along this
altitude and apply
case~\cC\ to each of the resulting triangles.

\cE\
Suppose $\D$ contains a vertex (say $B$) of $T$. If $\geo{AC}$
meets $\D$ then apply the previous case. So suppose $\geo{AC}$
misses $\D$ and consider the set of points $P$ of $\geo{BC}$
that are not joined to $A$ by a geodesic that misses $\D$. This
set is closed by lemma~\tag{thm-br-cover-basics}\rom5 and
nonempty because it contains $B$; we
let $B'$ be a point of this set closest to $C$. Since $B'\neq C$
there is a sequence of points in the interior of $\geo{B'C}$
approaching $B'$, each of which is joined with $A$ by a geodesic
missing $\D$. By lemma~\tag{thm-geodesics-converge}, there is a
geodesic $\geo{AB'}$. By the construction of $B'$, $\geo{AB'}$ meets $\D$.
Subdivide $T$ along this altitude and
apply  case \cB\ to $\tri{AB'C}$ and case \cD\ to $\tri{ABB'}$,
which of course reduces in turn to two applications of case \cC.

\cF\
Suppose $\D$ contains a point of $T$. By
lemma~\tag{thm-br-cover-basics}\rom3 there is a geodesic joining
this point with a vertex opposite it. Subdivide $T$ along this
altitude and apply case \cE\ to each of the resulting
triangles.

\cG\
Suppose an altitude of $T$ meets $\D$. Subdivide $T$ along
this altitude and apply case \cF\ to each of the resulting
triangles. 
\endproof

\beginproof{Proof of theorem~\tag{thm-br-cover-curvature}\rom2:}
Suppose that $\D$ is complete; we must show that $\Yhat$ is
geodesic. In light of lemma~\tag{thm-br-cover-basics}\rom3 it
suffices to show that any two points $x,z$ of $Y$ are joined
by a geodesic. We write $D$ for $d(x,z)$. Suppose first that
there exists a sequence  $y_i$ in $\D$ such that the sequence
$d(x,y_i)+d(y_i,z)$ converges to $D$. Then for each $i$ there is
a geodesic $\a_i$ (resp. $\b_i$) from $x$ (resp. $z$) to $y_i$
and 
we write $a_i$ (resp. $b_i$) for its length. By passing to a
subsequence we may suppose that the $a_i$ converge to a limit,
say $a$. Then the $b_i$ converge to $b=D-a$. Since neither $x$
nor $z$ lies in the closed set $\D$, the $a_i$ and $b_i$ are bounded away from
$0$, so $a>0$ and $b>0$. We will show that the $y_i$ form a
Cauchy sequence.

We let $\d>0$ be sufficiently small, by which we mean that
$\d<a$, $\d<b$, $a+\d<\pi/2\sx$ and $b+\d<\pi/2\sx$. This is
possible because each $a_i$ and $b_i$, hence each of $a$ and
$b$, is bounded above by $\diam(\Xhat)<\pi/2\sx$.  For such
$\d$, let $Y_\d=\set{y_i}{|a_i-a|<\d\hbox{ and }
|b_i-b|<\d}$. Now suppose $y_i,y_j\in Y_\d$ and let
$\c$ be the geodesic joining them. Since $\c$ lies in $\D$, the
geodesic triangle formed by $\a_i$, $\a_j$ and $\c$ satisfies
\cat, by lemma~\tag{thm-if-one-edge-in-branch-locus}. Similarly, the
triangle formed by $\b_i$, $\b_j$ and $\c$ also satisfies
\cat. We will show that if $\c$ were very long then its midpoint
would be problematic. Let $A$ be a closed annulus in $M_\x^2$
with inner radius $a-\d$ and outer radius $a+\d$. The conditions
we have imposed on $\d$ guarantee that the inner radius is
positive and (if $\x>0$) that $A$ lies in an open hemisphere. Let
$f(\d)$ be the maximum of the lengths of geodesic segments of
$M_\x^2$ that lie entirely in $A$. We observe that $f(\d)$ tends
to $0$ as $\d$ does.  Also, given any geodesic of $M_\x^2$ with
length $>2f(\d)$ and endpoints in $A$, its midpoint does not lie
in $A$ and hence lies at distance $<a$ from the center of
$A$.  (By
`the' center of $A$ when $\x>0$ we mean the center closer to
$A$.)  We define $g(\d)$ in a similar way, with $b$ in place of
$a$.  Now, if $\c$ were longer than $\max(2f(\d),2g(\d))$, then
by the \cat\ inequalities its midpoint $y'$ would satisfy
$d(x,y')<a$ and $d(y',z)<b$, a contradiction of the fact that
$d(x,z)=a+b$. 

We have shown that any two elements of $Y_\d$ have distance
bounded by $\max(2f(\d),2g(\d))$. Since this bound tends to $0$
as $\d$ does, the $y_i$ form a Cauchy sequence and hence
converge. Since the limit lies in $\D$, there are
geodesics joining it with
$x$ and with $z$. Concatenating these yields the required
geodesic.

Now suppose that there is no such sequence $y_i$. Then there is
a positive number $k$ such that $d(x,y)+d(y,z)>D+k$ for all $y$ in
$\D$, and there is also a sequence of paths in $Y$ from $x$ to
$z$ with lengths tending to $D$. From the first of these facts
we deduce that no path from $x$ to $z$ of length $<D+k$ meets
$\D$. So let $\c:[0,1]\to Y$ be a path from $x$ to $y$ of length
$<D+k$ and let $U$ be the set of $t\in[0,1]$ for which there is
not only a geodesic of $\Yhat$ from $x$ to $\c(t)$ but even one
that misses $\D$. The theorem follows because $U$ contains $0$
and is open (lemma~\tag{thm-br-cover-basics}\rom5) and
closed. To see that it is closed, let a sequence $t_i$ in $U$
converge to a point $t$ of $[0,1]$, and let $\b_i$ be a geodesic
from $x$ to $\c(t_i)$ that misses $\D$. By
lemma~\tag{thm-geodesics-converge}, the $\b_i$ converge to a
geodesic $\b$ from $x$ to $\c(t)$. The concatenation of $\b$ and
$\c|_{[t,1]}$ has length bounded by that of $\c$, which is less
than $D+k$. It follows that the concatenation does not meet
$\D$. In particular, $\b$ misses $\D$ and hence $t\in U$.
\endproof

\remark{Example:}
In some circumstances one can show that a branched cover is
\cat, even when the branch locus
is not complete. To illustrate this, we take $\Xhat$ to be the open unit
ball in $\R^3$ and $\D$ to be a diameter. With $Y$ as the
universal cover of $X=\Xhat-\D$, we can deduce that
$\Yhat=Y\cup\D$ is geodesic, even though $\D$ is not
complete. One simply takes $\bar\D$ to be the metric completion
of $\D$, $\bar X=X\cup\bar\D$, and considers the the branched
cover of $\bar X$ over $\bar\D$, where $X=\bar X-\bar\D$ is the
same as before, $Y$ is the universal cover of $X$, and $\bar
Y=Y\cup\bar\D$. 
That is, we add the endpoints of the diameter, then remove them
along with the diameter, take the cover as before, and then glue
the diameter and its endpoints back in.
Theorem~\tag{thm-br-cover-curvature} shows that
$\bar Y$ is CAT(0). Then, as an open ball in the CAT(0) space
$\bar Y$, $\Yhat$
is convex and hence also CAT(0).

We now present a local form of theorem~\tag{thm-br-cover-curvature}. The main
new feature is that the
projection map $\Yhat\to\Xhat$ is no longer required to be 1-1
on the branch locus. It also allows us to dispense with the
explicit diameter bounds for $\Xhat$ and $\Yhat$. We say that a metric space $\D$ is locally complete if
each of its points has a neighborhood whose closure is
metrically complete. This is equivalent to $\D$ being an open
subset of its completion. The pathological local properties of
the second example following
theorem~\tag{thm-br-cover-curvature} stem from the fact that the
branch locus used there is not locally complete.
If $\Xhat$ is a metric space and $\D\sset\Xhat$ then we say that
$\D$ is locally convex (in $\Xhat$) if each point has a
neighborhood $V$ such that $V\cap\D$ is convex in $\Xhat$.

\beginproclaim Theorem 
{\Tag{thm-local-br-cover}}. 
Suppose $\Xhat$ is a metric space of curvature $\leq\x$ for some
$\x\in\R$, and that $\D$ is a locally convex, locally
complete subset of $\Xhat$. Suppose $\Yhat$ is a metric space
and that $\pi:\Yhat\to\Xhat$ has the following
properties. First, each element of $\Dtil=\pi^{-1}(\D)$ has a
neighborhood $V$ such that $\pi|_V$ is a simple branched cover
of its image $\pi(V)$, over $\pi(V)\cap\D$. Second, $\pi$ is a local
isometry on $Y=\Yhat-\Dtil$. Then $\Yhat$ has curvature $\leq\x$.
\endproof
\endproclaim 

We omit the proof because we do not need the result and the
argument is a straightforward application of
theorem~\tag{thm-br-cover-curvature} and the idea used in the
example above.

\section{\Tag{sec-iterated}. 
Iterated branched covers of Riemannian manifolds}

In this section we define precisely what we mean by a branched
cover which is locally an iterated branched cover of a manifold
over a family of mutually orthogonal totally geodesic submanifolds. Then
we show that such a branched cover satisfies the
same upper bounds on local curvature as the base manifold. We
prove this only in the case of nonpositive curvature, but we
indicate what else is needed in the general case.

We say that a collection $\cals_0$ of subspaces of a real vector space $A$
is normal if the intersection of any $k\geq1$ members of $\cals_0$
has codimension $2k$. This means that each subspace has
codimension 2 and that they are as transverse as possible to
each other. The basic example is a subset of the coordinate
hyperplanes in $\C^n$, with $A$ being the underlying real vector
space. This is essentially the only example, in the following
sense. If $S_1,\ldots,S_n$ are the elements of $\cals_0$ then we
may introduce a basis $w_1,\ldots,w_m,x_1,y_1,\ldots,x_n,y_n$ of
$A$ such that each $S_i$ is the span of $w_1,\ldots,w_m$ and those
$x_j$ and $y_j$ with $j\neq i$. We will write $\S$ for the union
of the elements of $\S_0$.

Now suppose $\H_0$ is a family of immersed submanifolds of a
Riemannian manifold $\Mhat$ with union $\H$. We say that $\H_0$ is normal at
$x\in\Mhat$ if there is a family $\S_0$ of orthogonal
subspaces of $T_x\Mhat$ that are normal in the sense above and
have the following property. We require that there be an open
ball $U$ about $0$ in $T_x\Mhat$ which the exponential map $\exp_x$
carries diffeomorphically onto its image $V$, such that
$V\cap\H=\exp_x(U\cap\S)$, and such that $\exp(S\cap U)$ is a convex
subset of $V$ for each $S\in\S_0$. 
We say that $\H_0$ is normal if it is normal at
each $x\in\Mhat$. In this case, each element of $\H_0$ is
totally geodesic, distinct elements of $\H_0$ meet orthogonally
everywhere along their intersection,
and each self-intersection of an element of $\H_0$ is also
orthogonal. 

Let $x\in\Mhat$, $U$, $V$, $\S$ and $\H$ be as above, and write
$S_1,\ldots,S_n$ for the elements of $\S_0$. Then
$$ 
\pi_1(V-\H)\isomorphism\pi_1(U-\S)\isomorphism
\pi_1(T_x\Mhat-\S)\isomorphism \Z^n\;; 
$$
the first two isomorphisms are obvious and canonical, and the
last follows from the explicit description of $\S_0$ given
above. That is, $T_x\Mhat-\S$ is a product of $n$ punctured
planes and a Euclidean space.
We choose generators $\s_1,\ldots,\s_n$ for
$\pi_1(T_x\Mhat-\S)$ by taking a representative for $\s_i$ to be
a simple circular loop that links $S_i$ but none of the other
$S_j$. We say that a connected covering space of $T_x\Mhat-\S$
is standard if the subgroup of $\Z^n$ to which it
corresponds is generated by $\s_1^{d_1},\ldots,\s_n^{d_n}$ for some
$d_1,\ldots,d_n\in\Z$. We apply the same
terminology to the corresponding cover of $V-\H$.
In particular, the universal cover is
standard. An arbitrary covering space of $V-\H$ is called
standard if each of its components is. 

Now suppose $\Mhat$ is a Riemannian manifold and $\H_0$ is a
normal family of immersed submanifolds. We write
$M$ for $\Mhat-\H$. If $\pi:N\to M$ is a covering space then we
say that $N$ is a standard cover of $M$ if for each $x\in\Mhat$
with $V$ as above, $\pi:\pi^{-1}(V-\H)\to V-\H$ is a standard
covering in the above sense. In this case, we take $\Nhat$ to be
a certain subset of the metric completion of $N$: those points
of the completion which map to points of $\Mhat$ under the
completion of $\pi$. In particular, if $\Mhat$ is complete then
$\Nhat$ is the completion of $N$. We write $\pihat$ for the
natural extension $\Nhat\to\Mhat$ of $\pi$, and we say that this
map is a standard branched covering of $\Mhat$ over $\H_0$.

The simplest example of a standard branched cover is when
$\Mhat=\C^n$, $\H_0$ is the set of coordinate hyperplanes, $M$
is their complement and $\pi:N\to M$ is the covering space with
$N=\C^n-\H$ and
$\pi:(z_1,\ldots,z_n)\mapsto(z_1^{d_1},\ldots,z_n^{d_n})$. The
generalization to the case of locally infinite branching
requires the more complicated discussion in terms of metric
completions of covering spaces.
It is also possible that different components of the preimage of
$V-\H$ are inequivalent covering spaces of $V-\H$. This can
happen when $N\to M$ is an irregular cover.

We will need the following two general lemmas, whose proofs
should be skipped on a first reading. The first simplifies the
task of establishing local curvature conditions and the second
says that one may often ignore the added points when taking a
metric completion of a length space.

\beginproclaim Lemma 
{\Tag{thm-boring-lemma}}. 
Let $X$ be a length space with metric $d$. Let $Y$ be a
path-connected subset of $X$ with the property that
$$ 
\d(y,z)=\inf\set{\ell(\c)}{\hbox{\rm $\c$ a path in $Y$ joining $y$ and $z$}} 
$$
is finite for all $y,z\in Y$, so that $(Y,\d)$ is a length
space. Suppose also that $(Y,\d)$ is \cat\ for some
$\x\in\R$. Then any point of the interior of $Y$ admits a
neighborhood which is convex in $X$ and also \cat.
\endproclaim

\beginproof{Proof:}
We will define three open balls; all are balls with respect to
the metric $d$, rather than $\d$. Suppose $x$ lies in the
interior of $Y$ and that $U$ is an open ball centered at $x$ and
lying in $Y$. Let $r$ be the radius of $U$, and let $V$ be the
open ball with center $x$ and radius $r/2$. A simple argument
shows that $d(y,z)=\d(y,z)$ for all $y,z\in V$. Let $W$ be the
open ball with center $x$ and radius $r'=\min(r/4,\pi/\sx)$. Any
two points $y$, $z$ of $W$ are joined by a path $\c$ in $Y$ that
is a geodesic with respect to $\d$. Now, $\c$ lies in $V$ by the
triangle inequality, so $\c$ is also a geodesic with respect to
$d$, and $d(x,t)=\d(x,t)$ for all points $t$ of $\c$. By the
\cat\ inequality in $(Y,\d)$, applied to a triangle obtained by
joining $y$ and $z$ to $x$ with geodesics, we see that
$\d(x,t)<r'$ for all $t$. This shows that $\c$
lies in $W$. Since the same argument applies to every geodesic of
$X$ joining $y$ and $z$, we see that $W$ is convex in $X$. Since
$d$ and $\d$ coincide on $W$, $W$ is \cat.
\endproof

We say that the interior of a path $\c$ with domain $[a,b]$ lies
in a subset $Z$ of $\Xhat$ if $\c((a,b))\sset Z$.

\beginproclaim Lemma 
{\Tag{thm-paths-with-good-interiors}}.  Let $X$ be a length
space with metric $d$ and let $\Xhat$ be its metric
completion. For any $x,y\in\Xhat$ there are paths joining $x$
and $y$ with interiors in $X$ and lengths arbitrarily close to
$d(x,y)$. Furthermore, the intersection of $X$ with any open ball in
$\Xhat$ is path-connected.
\endproclaim

\beginproof{Proof:} 
Choose a sequence of points $x_i\in X$ that tend to $x$, such
that $\sum d(x_i,x_{i+1})<\infty$. By choosing short paths in
$X$ joining each $x_i$ to $x_{i+1}$ and concatenating them, we
obtain an open path $\c:[0,1)\to X$ of finite length, which can be
extended to $[0,1]$ by defining $\c(1)=x$. The extension is
continuous because $\c$ has finite length. By taking subpaths of
$\c$ we see that for all $\e>0$ there is a path of length $<\e$
from some point $x'$ of $X$ to $x$, with interior in $X$. The same
result holds with $y$ in place of $x$. Given $\e>0$, choose such
paths of lengths $<\e/4$. Then $x'$ and $y'$ may be
joined by a path in $X$ of length $<d(x,y)+2\e/4$. Putting the three
paths together establishes the first claim.

Now suppose $U\sset\Xhat$ is the open ball of radius $r$ and
center $x$, and that $y,z\in U\cap X$. By the above, there are paths from
$y$ and $z$ to $x$ with lengths $<r$ and interiors in $X$. These
paths obviously lie in $U$. If $y'$ and $z'$
are points in $U\cap X$ on these paths at distance $<r/2$ from
$x$, then they may be joined by a path in $X$ of length
$<r$. Such a path lies in $U$ by the triangle inequality. We
have joined $y$ to $y'$, $y'$ to $z'$ and $z'$ to $z$ by paths
in $U\cap X$, establishing the second claim.
\endproof

\beginproclaim Theorem 
\Tag{thm-iterated-br-cover}.  
If a Riemannian manifold $\Mhat$
has sectional curvature bounded above by $\x\leq0$ and
$\pihat:\Nhat\to\Mhat$ is a standard branched cover over a
normal family $\H_0$ of immersed submanifolds of $\Mhat$, then
$\Nhat$ is locally \cat.
\endproclaim 

\remark{Remark:}
For a global version of this result see theorem~\tag{thm-hyperplane-complements-aspherical}.

\beginproof{Proof:}
We will write $\Htil$ for $\pihat^{-1}(\H)$. Let
$\xtil\in\Nhat$, $x=\pihat(\xtil)$, and let $U$, $V$, $\S_0$ and $\S$ be
as in the definition of the normality of $\H_0$ at $x$.  We may
suppose without loss of generality that $\S_0\neq\emptyset$.
Let $r$ be the common radius of $U$ and $V$. Without loss of
generality we may take $r$ small enough so that $V$ and all
smaller balls centered at $x$ are convex in $\Mhat$. We
write $S_1,\ldots,S_n$ for the elements of $\S_0$, and $T_i$ for
$\exp_x(U\cap S_i)\sset V$. We choose $0<r'<r$ such that
the orthogonal projection maps from the closed $r'$-ball $B$ about
$x$ to the $T_i$ are well-behaved. By this we mean that for each
$i$, there is a fiberwise starshaped (about $0$) set in the restriction to
$T_i\cap B$ of the normal bundle of $T_i$, which is carried
diffeomorphically onto $B$ by the exponential map. The
orthogonal projection maps $B\to B\cap T_i$ are then obtained by
applying the inverse of this diffeomorphism followed by the
natural projection of the normal bundle to $T_i$. These maps
will not be used until late in the proof.

For $t\in[0,1]$ and $p\in V$ let $t.p$ denote the point in $V$
on the radial segment from $x$ to $p$ at distance $t\cdot
d(x,p)$ from $x$. Then the radial homotopy $\Gamma:[0,1]\times
V\to V$ given by $(t,p)\mapsto(1-t).p$ is a deformation
retraction of $V$ to $\{x\}$. Observe that if $p\in V-\H$ then
$t.p$ also lies in $V-\H$ for all $t\neq0$.  We may therefore 
lift $\Gamma|_{[0,1)\times(V-\H)}$ to an `open homotopy'
$[0,1)\times\pi^{-1}(V-\H)\to\pi^{-1}(V-\H)$ in the obvious
way. This is a lipschitz map and therefore extends  to a
homotopy $\Gammatil:[0,1]\times\pihat^{-1}(V)\to\pihat^{-1}(V)$. We
call this the radial homotopy; its tracks are geodesics of
$\Nhat$ and project to radial segments of $V$.

The first consequence of this analysis is that $\pihat^{-1}(V)$
is the union of the open $r$-balls about the points of
$\pihat^{-1}(x)$. The second consequence is that distinct
preimages of $x$ lies at distance $\geq 2r$ from each other. For
a path of length $<2r$ between  two preimages would lie in
$\pihat^{-1}(V)$ and then the deformation retraction of
$\pihat^{-1}(V)$ to $\pihat^{-1}(x)$ shows that the endpoints of
the path coincide. This implies that the open $r$-ball
$\Vtil$ about $\xtil$ is a component of $\pihat^{-1}(V)$, and
therefore the restriction of $\pi$ to $\Vtil-\Htil$ is a
covering map. We note  that $\Vtil-\Htil$ is connected, by
lemma~\tag{thm-paths-with-good-interiors}. The radial homotopy
also shows that the closed
$r'$-ball $\Btil$ about $\xtil$ is the preimage in $\Vtil$ of
$B$, that $\Btil$ is path-connected, and that for all
$\ytil,\ztil\in\Btil$,
$$ 
\d(\ytil,\ztil)=\inf\set{\ell(\c)}
{\hbox{$\c$ a path in $\Btil$ joining $\ytil$ and $\ztil$}} 
$$
is finite. Finally, since $\Vtil-\Htil$ is connected, the radial
homotopy shows that $\Btil-\Htil$ is also connected.

To show that $\xtil$ admits a convex \cat\ neighborhood it
suffices by lemma~\tag{thm-boring-lemma} to show that
$(\Btil,\d)$ is \cat. We will prove this by realizing $\Btil$ as
an iterated simple branched cover of $B$. We let
$\s_1,\ldots,\s_n$ denoted generators for
$G=\pi_1(B-\H)=\pi_1(T_x\Mhat-\S)\isomorphism\Z^n$ of the sort
discussed above. Since $\pi:N\to M$ is a standard covering,
there are $d_1,\ldots,d_n\in\Z$ such that the subgroup of $G$
associated to the covering 
$\Btil-\Htil\to B-\H$ is generated by
$\s_1^{d_1},\ldots,\s_n^{d_n}$. For each $k=0,\ldots,n$, let
$G_k$ be the subgroup generated by
$\s_1^{d_1},\ldots,\s_k^{d_k},\s_{k+1},\ldots,\s_n$. We let $B_k$ be
the metric completion of the cover of $B-\H$ associated to
$G_k$, equipped with the natural path
metric. Then $B_k$ is the standard branched cover of $B$,
branched over the $T_i\cap B$, with branching indices
$d_1,\ldots,d_k,1,\ldots,1$. In particular, $B_0=B$ and
$B_n=(\Btil,\d)$.

We write $p_k$ for the natural projection $B_k\to B$
obtained by extending the covering map to a map of metric completions. 
Because
$G_{k+1}\sset G_k$, there is a covering map 
$B_{k+1}-p_{k+1}^{-1}(\H)\to 
B_{k}-p_{k}^{-1}(\H)$ whose completion
$q_{k+1}:B_{k+1}\to B_k$ satisfies $p_k\circ
q_{k+1}=p_{k+1}$. 
For each
$k=0,\ldots,n-1$ we let $\D_k=p_k^{-1}(T_{k+1})$. 
We will show that
$q_{k+1}$ is a simple branched covering map with branch locus
$\D_k$, for each $k=0,\ldots,n-1$. 

First we claim that $q_{k+1}$ carries $q_{k+1}^{-1}(\D_k)$
bijectively to $\D_k$ and is a covering map on the complement
of $q_{k+1}^{-1}(\D_k)$. To see this, observe that $B$ is
bilipshitz to the metric product $A_1\times\cdots\times A_n\times D$ of $n$
closed Euclidean disks and a closed Euclidean ball, such that
$T_i\cap B$ is identified with the set of points in the product
whose $i$th coordinate is  the
center (say $0$) of $A_i$. This identifies each $B_k$ with 
$\Atil_1\times\cdots\times\Atil_k\times
A_{k+1}\times\cdots\times A_n\times D$, where $\Atil_j$ is the
metric completion of the $d_j$-fold cover of $A_k-\{0\}$ (or the
universal cover if $d_j=0$).
The metric completion of this cover of $A_j-\{0\}$ is obtained
from the cover by adjoining a single point, which lies over $0$.
Then $q_{k+1}:B_{k+1}\to B_k$ is
given by the branched cover $\Atil_{k+1}\to
A_{k+1}$ and the identity maps on $\Atil_1,\ldots,\Atil_{k}$,
$A_{k+2},\ldots,A_n$, and $D$. The claim is now obvious.

Now we show that $q_{k+1}$ is a local isometry away from
$q_{k+1}^{-1}(\D_k)$. If
$\ytil\in B_{k+1}-q_{k+1}^{-1}(\D_k)$ then by the previous
paragraph there is an $s>0$ and a neighborhood $E$ of $\ytil$
such that $q_{k+1}$ carries $E$ homeomorphically onto its image,
which is the $s$-ball about $q_{k+1}(\ytil)$. It is then easy to
see that $q_{k+1}$ carries the open $(s/2)$-ball about $\ytil$
isometrically onto its image.

To complete the proof that $q_{k+1}$ is a simple branched cover,
we need only show that $d(\ytil,\ztil)=d(y,z)$, when
$\ytil,\ztil\in B_{k+1}$ have images $y$ and $z$ under
$q_{k+1}$ and at least one of $y$ and $z$ lies in $\D_k$. Without
loss of generality we may take $z\in\D_k$, and by continuity it
suffices to treat the case $y\notin\D_k$. It is obvious that
$d(y,z)\leq d(\ytil,\ztil)$. To see the converse, let $\c_i$ be
a sequence of paths in $B_k$ from $y$ to $z$, with interiors in
$B_k-\D_k$ and lengths approaching $d(y,z)$. This is possible by
lemma~\tag{thm-paths-with-good-interiors}. By lifting each path
except for its final endpoint,
and extending using completeness, we obtain paths $\ctil_i$ from
$\ytil$ to points of $B_{k+1}$ lying over $z$, with the same
lengths as the $\c_i$. By the injectivity of $q_{k+1}$ on
$q_{k+1}^{-1}(\D_k)$, all of these are paths from $\ytil$ to
$\ztil$, establishing the claim.

In order to use theorem~\tag{thm-br-cover-curvature}
inductively, we will need to know that $\D_k$ is a convex subset of
$B_k$. We will use the orthogonal projections
introduced earlier. Each of these projections $B\to B\cap T_j$
may be realized by a deformation retraction along geodesics. The
retraction is distance non-increasing, since each $T_j$ is
totally geodesic and $\Mhat$ has sectional curvature
$\leq0$. Because $T_1,\ldots,T_{k-1}$ are totally geodesic and
orthogonal to $T_k$, the track of the deformation retraction to
$T_k$ starting at a point outside $T_1\cup\cdots\cup T_{k-1}$
misses $T_1\cup\cdots\cup T_{k-1}$ entirely. Therefore the
deformation lifts to a deformation retraction of
$B_k-p_k^{-1}(T_1\cup\cdots\cup T_{k-1})$ to
$\D_k-p_k^{-1}(T_1\cup\cdots\cup T_{k-1})$. This extends to a
distance nonincreasing retraction $B_k\to\D_k$, which we will
also call an orthogonal projection.

Now we prove by simultaneous induction that $B_k$ is \cat\
and that $\D_k$ is convex in $B_k$. The fact that
$B_0=B$ is \cat\ follows from its convexity in $V$ and the
fact that $V$ is \cat, which in turn follows from the proof of
\ecite{troyanov:spaces-of-negative-curvature}{theorem~12}. That
theorem asserts that simply connected complete Riemannian
manifolds of sectional curvature $\leq\x$ are \cat, but the
proof shows that the completeness condition may be replaced by
the weaker condition that the manifold be geodesic. The convexity of
$\D_0=T_1\cap B$ in $B$ follows from the convexity of $T_1$ in
$V$. Now the inductive step is easy. If $B_k$ is \cat\ and
$\D_k$ is convex in $B_k$ then $B_{k+1}$ is \cat\ by
theorem~\tag{thm-br-cover-curvature}. 
In particular, geodesics in $B_{k+1}$ are unique.
Then if $\c$ is a
geodesic of $B_{k+1}$ with endpoints in $\D_{k+1}$,  the
orthogonal projection to $\D_{k+1}$ carries $\c$ to a path of
length $\leq\ell(\c)$ with the same endpoints. By the
uniqueness of geodesics, $\c$ lies in
$\D_{k+1}$, so we have proven that $\D_{k+1}$ is convex in
$B_{k+1}$. The theorem follows by induction.
\endproof

\remark{Remark:}
We indicate here the additional work required to prove the theorem
when $\x>0$. A minor point is that one should take $r<\pi/4\sx$,
so that all of the $B_k$ have diameters $<\pi/2\sx$. A more
substantial change is required where we used the fact that the
orthogonal projection maps $B\to B\cap T_j$ are distance
non-increasing, because this fails in the presence of positive curvature. All that is important here is that the length of a path
in $B$ {\it with endpoints in $T_k$} does not increase under
projection to $T_k$. Even this is not true, but we only
need the result for paths of length $<2r$. One should choose
$r'$ small enough so that any path in $B$ of length $<2r$ with
endpoints in $T_k$ grows no longer under the projection to
$T_k$. Presumably this can be done but I have not checked the
details. 

Theorem~\tag{thm-iterated-br-cover} has been widely believed,
but  ours seems to be the first proof.  The theorem overlaps
partly with theorem~5.3 of
\cite{charney-davis:branched-covers-of-riemannian-manifolds},
which considers locally finite branched covers of Riemannian
manifolds over subsets considerably more complicated than
mutually orthogonal submanifolds. Unfortunately there is a gap
in the proof of that theorem which I do not know how to bridge (lemma~5.7 does
not seem to follow from lemma~5.6). Nevertheless I regard
the `infinitesimal'
\cat\ condition (condition~3 of  theorem~5.3) 
as very natural, and expect that the theorem not only
holds but extends to the case of locally infinite branching.

\section{\Tag{sec-apps}. Applications}

In this section we solve the problems which motivated our
investigation, concerning the asphericity of certain moduli
spaces. By using known models for the moduli spaces of cubic
surfaces in $\cp^3$ and Enriques surfaces we will show that both
of these spaces have contractible universal covers. In both
cases a key ingredient is the following theorem, which is a sort
of global version of theorem~\tag{thm-iterated-br-cover}.

\beginproclaim Theorem 
{\Tag{thm-hyperplane-complements-aspherical}}. 
Let $\Mhat$ be a complete simply connected Riemannian manifold
with section curvature bounded above by $\x\leq0$. Let $\H_0$ be
a family of complete submanifolds which are normal in the sense
of section~\tag{sec-iterated}, and let $\H$ be the union of the members of
$\H_0$. Then the metric completion $\Nhat$ of the universal
cover $N$ of $\Mhat-\H$ is \cat, and $N$ and $\Nhat$ are
contractible. 
\endproclaim 

\beginproclaim Corollary 
{\Tag{thm-cubic-moduli-space}}. 
The moduli space $\moduli$ of smooth cubic surfaces in $\cp^3$ is
aspherical. 
\endproclaim 

\beginproof{Proof:}
We begin by recalling the main result of \cite{allcock:ch4-cubic-moduli}. We let $\w$ be
a primitive cube root of unity and set $\E=\Z[\w]$, a discrete
subring of $\C$. Let $L$ be the lattice $\E^5$ equipped with the
Hermitian inner product
$$ 
h(x,y)=-x_0\bar y_0+x_1\bar y_1+\cdots+x_4\bar y_4\quad. 
$$
Then the complex hyperbolic space 
$\ch^4$ may be taken to be the set of lines in $\C^5$ on
which $h$ is negative-definite, so that $\ch^4$ is a subset of
$\cp^4$. Let $\H_0$ be the set of (complex) hyperplanes in $\ch^4$
which are the orthogonal complements of those $r\in L$ with
$h(r,r)=1$. Let $\Gamma$ be the unitary group of $L$, which
is obviously discrete in $\U(4,1)$. By \cite{allcock:ch4-cubic-moduli},
$(\ch^4-\H)/\Gamma$ is isomorphic as an orbifold to $\moduli$.  

To see that $\H_0$ is locally finite, observe that $\Gamma$
contains a complex reflection in each element $H$ of $\H_0$,
That is, if $H$ corresponds to $r\in L$ with $h(r,r)=1$, then
$\Gamma$ contains an element fixing $r^\perp$ (and hence $H$)
pointwise and multiplying $r$ by $-1$. 
If $\H_0$ failed to be locally finite then the existence
of these reflections would contradict the discreteness of
$\Gamma$. Now consider two elements of
$\H_0$ that meet in $\ch^4$, and suppose that they are
associated to $r,r'\in L$. Since they meet, $h$ is
positive-definite on the span of $r$ and $r'$. Since
$h(r,r)=h(r',r')=1$, positive-definiteness requires
$|h(r,r')|<1$. Since $h(r,r')\in\E$ we must have
$h(r,r')=0$. This shows that any two elements of $\H_0$ that
meet do so orthogonally. 

Since $\ch^4$ has negative sectional curvature,
theorem~\tag{thm-hyperplane-complements-aspherical} implies that
$\ch^4-\H$ has contractible universal cover. This is also the
orbifold universal cover of $(\ch^4-\H)/\Gamma$, so the result
follows.  
\endproof

\beginproclaim Corollary 
{\Tag{thm-enriques-moduli-space}}. 
The period space for smooth complex Enriques surfaces (defined
below) is aspherical.
\endproclaim 

\beginproof{Proof:}
This is similar to the previous proof. By the Torelli theorem
for Enriques surfaces (\cite{horikawa:periods-of-enriques-surfaces-I} and \cite{namikawa:periods-of-enriques-surfaces}), the isomorphism
classes of smooth complex Enriques surfaces are in 1-1
correspondence with the points of the period space
$(\DD-\H)/\Gamma$. Here $\DD$ is the symmetric space for
$O(2,10)$, $\Gamma$ is a certain discrete subgroup, and $\H_0$
is a certain $\Gamma$-invariant arrangement of complex
hyperplanes. By \cite{allcock:period-lattice-for-enriques-surfaces}, $\Gamma$ may be taken to be the
isometry group of the lattice $L$ which is $\Z^{12}$ equipped
with the inner product
$$ 
x\cdot y=x_1y_1+x_2y_2-x_3y_3-\cdots-x_{12}y_{12}\;. 
$$
A concrete model for $\DD$ is the set of $v\in L\tensor\C$
satisfying $v\cdot v=0$ and $v\cdot\bar{v}>0$. $\H_0$ may be
taken to be the set of (complex) hyperplanes in $\DD$ which are
the orthogonal complements of the norm $-1$ vectors of
$L$. The arguments of the previous proof show that
$\H_0$ is normal. As a symmetric space of noncompact type, $\DD$
has sectional curvature~$\leq0$, so that theorem~\tag{thm-hyperplane-complements-aspherical}
applies and $\DD-\H_0$ is aspherical. It follows that the period
space, the orbifold $(\DD-\H)/\Gamma$, is also aspherical.
\endproof

\remark{Remark:}
We have referred to the period space of Enriques surfaces rather
than to a moduli space. This is because it is highly nontrivial to
assemble the isomorphism classes of Enriques surfaces into a
moduli space $\moduli$. One way to do this is to equip the
surfaces with suitable extra structure and then use geometric
invariant theory, as in \cite{shah:projective-degeneration-of-enriques-surfaces}. Then $\moduli$ has a natural
topology and is identified with some finite cover $C$ of
$(\DD-\H)/\Gamma$. Strictly speaking, one should also impose
additional structure to make sure that $C$ is a manifold
and not just an orbifold. The reason is that the orbifold
structure on $C$ is not terribly relevant to $\moduli$. (There 
are Enriques surfaces with infinitely many automorphisms,
as well as those with only finitely many, so there is not much
hope of a reasonable orbifold structure on any $\moduli$.) But
if $C$ is a manifold then the homeomorphism of $C$ with
$\moduli$ shows that $\moduli$ is aspherical.

Now we will prove theorem~\tag{thm-hyperplane-complements-aspherical}.  In the proofs we will conform
to the notation of section~\tag{sec-iterated} by writing $M$ for
$\Mhat-\H$, $\pi$
for the covering map $N\to M$, $\pihat$ for its completion, and
$\Htil$ for $\pihat^{-1}(\H)\sset\Nhat$.

\beginproclaim Lemma 
{\Tag{thm-neighborhoods-described}}.  The map $\pi$ is a
standard covering and the map $\pihat$ is a standard branched
cover over $\H$. More precisely, suppose $\xtil\in\Nhat$,
$x=\pihat(x)$, $V$ is an open ball about $x$ that meets no
element of $\H_0$ except those passing through $x$, and $\Vtil$
is the ball of the same radius about $\xtil$. Then $\Vtil-\Htil$
is a copy of the universal cover of $V-\Htil$.
\endproclaim

\beginproof{Proof:}
We prove only the last assertion. Let $\I_0$ be
the set of elements of $\H_0$ which contain $x$, and let $\I$ be
their union. Then there is a natural sequence of group
homomorphisms 
$$ 
\pi_1(V-\I)
\to\pi_1(\Mhat-\H)
\to\pi_1(\Mhat-\I)
\to\pi_1(V-\I)
\quad.
$$
The first and second maps are induced by inclusions of the
indicated spaces and the third is induced by a retraction of
$\Mhat-\I$ to $V-\I$ along geodesic rays based at $x$. The
composition is obviously the identity map, which shows that  
$\pi_1(V-\I)\to\pi_1(\Mhat-\H)$ is
injective. Therefore each component of the preimage in $N$ of
$V-\H$ is a copy of the universal cover of $V-\H$.  By
lemma~\tag{thm-paths-with-good-interiors}, $\Vtil-\Htil$ is
connected, and as in the proof of
theorem~\tag{thm-iterated-br-cover} it is a component of
$\pi^{-1}(V-\H)$. This completes the proof.
\endproof

\beginproclaim Lemma 
{\Tag{thm-weak-homotopy-equiv}}. 
The inclusion $N\to\Nhat$ is a weak homotopy equivalence.
\endproclaim

\beginproof{Proof:}
First we show that for each $\xtil\in\Nhat$ there is a homotopy
of $\Nhat$ to itself that carries $N$ into itself and also some
neighborhood of $\xtil$ into $N$. We write $n$ for the number of
hyperplanes passing through $x=\pi(\xtil)$. By the previous
lemma and the ideas of the proof of theorem~\tag{thm-iterated-br-cover}, there is a
closed neighborhood $\Vtil$ of $\xtil$ which is homeomorphic to
the metric product $(\Atil)^n\times D$, where $\Atil$ is the
metric completion of the universal cover of a closed Euclidean disk
minus its center and $D$ is a closed Euclidean ball. It is easy
to see that $\Atil$ is homeomorphic to a `wedge' in the
plane, by which we mean 
$$ 
\Atil\isomorphism
\{(0,0)\}\;\cup\;
\set{(x,y)\in\R^2}{0<x,\; |y|<x,\; x^2+y^2\leq1}
\quad.
$$
There is obviously a homotopy of $\Atil$ into $\Atil-\{(0,0)\}$
which is supported on a small neighborhood of $(0,0)$. We obtain
the desired homotopy of $\Nhat$ by applying this homotopy to
each factor $\Atil$ of $\Vtil$ and fixing each point of
$\Nhat-\Vtil$.

Now, if $f:S^k\to\Nhat$ represents any element of the homotopy
group $\pi_k(\Nhat)$
then we may cover $f(S^k)$ with finitely many open sets, each of
which is carried into $N$ by some homotopy of $\Nhat$ that also
carries all of $N$ into $N$. Applying these homotopies one after
another shows that $f$ is homotopic to a map $S^k\to N$. This
shows that $\pi_k(N)$ surjects onto $\pi_k(\Nhat)$ for all
$k$.  The same argument applied to balls rather than
spheres shows that $\pi_k(N)$ also injects, completing the
proof.
\endproof

\beginproof{Proof of theorem~\tag{thm-hyperplane-complements-aspherical}:}
By lemma~\tag{thm-neighborhoods-described}, $\Nhat$ is a
standard branched cover of $\Mhat$ over the normal family $\H_0$. Since $\Mhat$ has
sectional curvature $\leq\x\leq0$,
theorem~\tag{thm-iterated-br-cover} shows that $\Nhat$ is
locally \cat. Since $N$ is simply connected
lemma~\tag{thm-weak-homotopy-equiv} implies that $\Nhat$ is
also. Theorem~\tag{thm-cartan-hadamard} now implies that $\Nhat$
is \cat\ and hence contractible. In particular, all of its
homotopy groups vanish, and by another application of
lemma~\tag{thm-weak-homotopy-equiv} the same is true of $N$. As
a manifold all of whose homotopy groups vanish, $N$ is
contractible.
\endproof

In the introduction we promised to show that the inclusion
$N\to\Nhat$ is a homotopy equivalence, but all we have
established so far is a weak homotopy equivalence. The stronger
result follows immediately from
theorem~\tag{thm-hyperplane-complements-aspherical}, since any
inclusion of one contractible space into another is a homotopy
equivalence.

\section{References}


\bibitem{alexandrov51:theorem_on_triangles}
A.~D. Alexandrov.
 A theorem on triangles in a metric space and some of its
  applications.
 {\it Trudy Math. Inst. Steks.}, 38:5--23, 1951.

\bibitem{allcock:ch4-cubic-moduli}
D.~Allcock, J.~Carlson, and D.~Toledo.
 The complex hyperbolic geometry of the moduli space of cubic
  surfaces.
 Preprint, 1998.

\bibitem{allcock:period-lattice-for-enriques-surfaces}
D.~J. Allcock.
 The period lattice for {E}nriques surfaces.
 Preprint 1999.

\bibitem{bridson95:metric_spaces_of_nonpositive_curvature}
M.~Bridson and A.~Haefliger.
 Metric spaces of non-positive curvature.
 Book preprint, 1995.

\bibitem{charney-davis:branched-covers-of-riemannian-manifolds}
R.~Charney and M.~Davis.
 Singular metrics of nonpositive curvature on branched covers of
  {R}iemannian manifolds.
 {\it Am. J. Math.}, 115:929--1009, 1993.

\bibitem{davis91:hyperbolization_of_polyhedra}
M.~Davis and T.~Januszkiewicz.
 Hyperbolization of polyhedra.
 {\it J. Diff. Geom.}, 34:347--388, 1991.

\bibitem{gp-thy-geom-viewpoint}
E.~Ghys et~al., editors.
 {\it Group Theory from a Geometrical Viewpoint}. World Scientific,
  1991.

\bibitem{gromov:hyperbolic-groups}
M.~Gromov.
 Hyperbolic groups.
 In S.~M. Gersten, editor, {\it Essays in Group Theory}, volume~8 of
  {\it MSRI Publications}, pages 75--263. Springer-Verlag, 1987.

\bibitem{horikawa:periods-of-enriques-surfaces-I}
E.~Horikawa.
 On the periods of {E}nriques surfaces. {I}.
 {\it Math. Ann.}, 234:73--88, 1978.

\bibitem{januszkiewicz:hyperbolizations}
T.~Januszkiewicz.
 Hyperbolizations.
 In Ghys et~al. \cite{gp-thy-geom-viewpoint}, pages 464--490.

\bibitem{namikawa:periods-of-enriques-surfaces}
Y.~Namikawa.
 Periods of {E}nriques surfaces.
 {\it Math. Ann.}, 270:201--222, 1985.

\bibitem{paulin:hyperbolic_groups_via_hyperbolizing_polyhedra}
F.~Paulin.
 Construction of hyperbolic groups via hyperbolizations of polyhedra.
 In Ghys et~al. \cite{gp-thy-geom-viewpoint}, pages 313--372.

\bibitem{shah:projective-degeneration-of-enriques-surfaces}
J.~Shah.
 Projective degenerations of {E}nriques' surfaces.
 {\it Math. Ann.}, 256(4):475--495, 1981.

\bibitem{gromov-thurston:pinching-constants-for-hyperbolic-manifolds}
W.~Thurston and M.~Gromov.
 Pinching constants for hyperbolic manifolds.
 {\it Invent. Math.}, 89(1):1--12, 1987.

\bibitem{troyanov:spaces-of-negative-curvature}
M.~Troyanov.
 Espaces {\`a} courbure n{\'e}gative et groups hyperboliques.
 In E.~Ghys and P.~de~la Harpe, editors, {\it Sur les Groupes
  Hyperboliques d'apr{\`e}s {M}ikhael {G}romov}, volume~83 of {\it Progress in
  Mathematics}, pages 47--66. Birk{\"a}user, 1990.


\bigskip
\parindent=0pt
Department of Mathematics, Harvard University

One Oxford Street, Cambridge, MA 02138.

email: {\it allcock@math.harvard.edu}

web page: {\it http://www.math.harvard.edu/$\sim$allcock}

\bye